\documentclass[preprint, 11pt]{elsarticle}

\usepackage[T1]{fontenc}
\usepackage[utf8]{inputenc}
\usepackage[hmargin=0.9in]{geometry}
\usepackage[babel,german=quotes]{csquotes}
\usepackage{subfig}
\usepackage{graphicx}
\usepackage[colorlinks=true,citecolor=blue,urlcolor=blue]{hyperref}
\usepackage{caption}
\usepackage{amsmath}          
\usepackage{amsthm}           
\usepackage{amssymb}          
\usepackage{bm, amsfonts}     
\usepackage{xcolor}
\usepackage{fixmath}
\usepackage{tikz}
\usepackage{bm}
\usepackage{makecell}

\newtheorem{thm}{Theorem}
\newdefinition{example}{Example}
\newproof{pf}{Proof}
\newcommand{\proofend}{\hfill $\Box$}

\newcommand {\dx} {\,{\rm d}{\mathbf x}}

\newcommand {\ds} {\,{\rm d}{\mathrm s}}

  \newcommand{\R}{\mathbb{R}}
  \newdefinition{rmk}{Remark}
  
  \newcommand{\pd}[2]{\frac{\partial #1}{\partial #2}}

  \newcommand{\td}[2]{\frac{\mathrm d #1}{\mathrm d #2}}

\newcommand{\beq}{\begin{equation}}
\newcommand{\eeq}{\end{equation}}

\makeatletter
\def\ps@pprintTitle{%
  \let\@oddhead\@empty
  \let\@evenhead\@empty
  \def\@oddfoot{
    \footnotesize\itshape
    \hfill\today
  }%
  \let\@evenfoot\@oddfoot}
\makeatother

\begin{document}

\begin{frontmatter} 
  \title{Entropy stabilization and property-preserving limiters for discontinuous Galerkin discretizations of nonlinear hyperbolic equations}

\author{Dmitri Kuzmin}
\ead{kuzmin@math.uni-dortmund.de}

\address{Institute of Applied Mathematics (LS III), TU Dortmund University\\ Vogelpothsweg 87,
  D-44227 Dortmund, Germany}

\journal{Journal of Computational Physics}

\begin{abstract}
  The methodology proposed in this paper bridges the gap between entropy stable and positivity-preserving discontinuous Galerkin (DG) methods for nonlinear hyperbolic problems. The entropy stability property and, optionally, preservation of local bounds for the cell averages are enforced using flux limiters based on entropy conditions and discrete maximum principles, respectively. Entropy production by the (limited) gradients of the piecewise-linear DG approximation is constrained using Rusanov-type entropy viscosity, as proposed by Abgrall in the context of nodal finite element approximations. We cast his algebraic entropy fix into a form suitable for arbitrary polynomial bases and, in particular, for modal DG approaches. The Taylor basis representation of the entropy stabilization term reveals that it penalizes the solution gradients in a manner similar to slope limiting and requires semi-implicit treatment to  achieve the desired effect. The implicit Taylor basis version of the Rusanov entropy fix preserves the sparsity pattern of the element mass matrix. Hence, no linear systems need to be solved if the Taylor basis is orthogonal and an explicit treatment of the remaining terms is adopted. The optional application of a vertex-based slope limiter constrains the piecewise-linear DG solution to be bounded by local maxima and minima of the cell averages. The combination of entropy stabilization with flux and slope limiting leads to constrained approximations that possess all desired properties. Numerical studies of the new limiting techniques and entropy correction procedures are performed for two scalar two-dimensional
  test problems with nonlinear and nonconvex flux functions.

\end{abstract}
\end{frontmatter}

\section{Introduction}

In recent years, significant advances have been made in the analysis and
design of property-preserving high-resolution finite element schemes for
hyperbolic problems. The essential properties of a physics-compatible
approximation include entropy stability and positivity preservation.
Entropy stable discontinuous Galerkin (DG) methods \cite{chen,pazner} are
usually derived using entropy conservative numerical fluxes and
additional dissipative terms (fluxes and/or element contributions)
depending on the gradients of entropy variables. Tadmor's seminal
work \cite{tadmor87} provides a general entropy stability criterion
for the analysis and design of such schemes \cite{tadmor,tadmor2016}. In
the case of a piecewise-linear or higher-order DG approximation, excess
entropy production by the gradients of the conserved quantities must
be balanced using nonlinear artificial diffusion operators and/or
limiters. The entropy correction term proposed by Abgrall \cite{abgrall}
and its generalizations presented in \cite{ranocha} penalize the
deviations of entropy variables from their cell averages using Rusanov-type
entropy viscosity. Artificial diffusion operators of this kind are
also widely used to enforce local discrete maximum principles,
preservation of invariant domains, and/or positivity preservation
in low-order components of residual distribution methods
\cite{Abgrall2006,RD-BFCT} and algebraic flux correction (AFC) schemes
\cite{Guermond2014,convex,afclps} for hyperbolic problems. The
use of flux and slope limiters makes it possible to adjust the
amounts of artificial diffusion or the gradients of the
numerical solution in an adaptive manner. A variety of
algebraic \cite{DG-BFCT,badia,Guermond2019,RD-BFCT} and
geometric \cite{barthjesp,dglim} limiting techniques can be
found in the literature on finite volume and DG methods for
hyperbolic conservation laws. The most recent approaches
are backed by theoretical proofs of positivity preservation
for cell averages and/or solution values at certain control
points \cite{giuliani,Guermond2018,IRP-DG,Moe2017,zhang1,zhang2}.
However, positivity-preserving high-resolution schemes may
converge to wrong weak solutions if they are not entropy
stable \cite{Guermond2014,entropyCG}. Conversely, an entropy stable high-order method may
exhibit excellent convergence behavior in smooth regions but
produce undershoots and overshoots in shock regions.

An algebraic limiting framework that ensures both entropy
stability and preservation of local bounds was introduced
in \cite{entropyCG} in the context of AFC schemes for
continuous Galerkin methods. In the present paper,
we constrain piecewise-linear ($\mathbb{P}_1$) Taylor basis
DG discretizations using similar flux correction tools
in addition to slope limiting. The proposed flux limiter
guarantees that the cell averages satisfy a semi-discrete
entropy inequality and a local maximum principle. The rates
of entropy production and dissipation inside mesh cells
are balanced using built-in gradient penalization which
corresponds to a diagonal Taylor basis form of Abgrall's
\cite{abgrall} correction term. In the process of Runge-Kutta
time integration, we treat this term implicitly and apply a
vertex-based version \cite{dglim} of the Barth-Jespersen
\cite{barthjesp} slope limiter. The final DG-$\mathbb{P}_1$
approximation stays in the range determined by the local
maxima and minima of property-preserving cell averages.
In Sections \ref{sec:ES}-\ref{sec:SL}, we present the new correction tools
and explain the underlying design principles. The numerical
examples of Section \ref{sec:num} illustrate the implications of entropy
stability and the capability of the proposed algorithms to
enforce the desired properties. We close this paper with
a summary of the results and possible extensions in Section \ref{sec:end}.

\section{Entropy stabilization of DG schemes}
\label{sec:ES}

Let $u(\mathbf{x},t)$ be a scalar conserved
quantity depending on the space location $\mathbf{x}\in
\R^d,\ d\in\{1,2,3\}$ and time instant $t\ge 0$. Consider
an initial value problem of the form 
\begin{subequations}
\begin{align}
 \pd{u}{t}+\nabla\cdot\mathbf{f}(u)=0 &\qquad\mbox{in}\ \R^d\times\R_+,
\label{ibvp-pde}\\
 u(\cdot,0)=u_0 &\qquad\mbox{in}\ \R^d,\label{ibvp-ic}
\end{align}
\end{subequations}
where $\mathbf{f}=(\mathsf{f}_1,\ldots,\mathsf{f}_d)$ is a possibly
nonlinear flux function and $u_0:\R^d\to\mathcal G$ is an initial data
belonging to a convex set $\mathcal G$. The set $\mathcal G$ is
called an {\it invariant set} of problem  \eqref{ibvp-pde}--\eqref{ibvp-ic}
if the exact solution $u$ stays in $\mathcal G$ for all $t> 0$ \cite{Guermond2016}.
A convex function
$\eta:\mathcal G\to\R$ is called an {\it entropy} and $v=\eta'$ is called an
{\it entropy}
variable if there exists an {\it entropy flux} 
$\mathbf{q}:\mathcal G\to\mathbb{R}^d$ such that $v(u)\mathbf{f}'(u)=\mathbf{q}'(u)$.
 A weak solution
$u$ of \eqref{ibvp-pde} is called an {\it entropy solution}
if the entropy inequality
\beq
 \pd{\eta}{t}+\nabla\cdot\mathbf{q}(u)\le 0 \qquad\mbox{in}\ \R^d\times\R_+
\label{ent-ineq}
\eeq
holds for any {\it entropy pair} $(\eta,\mathbf{q})$. For any smooth weak
solution, the conservation law
\beq
 \pd{\eta}{t}+\nabla\cdot\mathbf{q}(u)=0 \qquad\mbox{in}\ \R^d\times\R_+
 \eeq
 can be derived from \eqref{ibvp-pde} using multiplication by the entropy
 variable $v$, the chain rule, and the definition of an entropy pair. Hence,
 entropy is conserved in smooth regions and dissipated at shocks.

 A numerical scheme is called  {\it positivity-preserving}
 \cite{Moe2017,zhang1,zhang2} or, more formally,
 {\it invariant domain preserving} (IDP)
 \cite{Guermond2018,Guermond2016,IRP-DG,convex}
     if the solution of the (semi-)discrete problem
 is guaranteed to stay in an invariant set $\mathcal G$. A~well-designed
 discretization of \eqref{ibvp-pde} should also be {\it entropy stable}, i.e., 
a (semi-) discrete version of the entropy inequality \eqref{ent-ineq}
should hold. The lack of entropy stability is a typical reason for
convergence of numerical methods to nonphysical weak solutions.

Let us discretize \eqref{ibvp-pde} in space using a piecewise-linear
DG approximation  $u_h$ on a computational mesh $\mathcal T_h$
consisting of $E_h$ elements. The restriction of $u_h$ to element 
$K_i,\ i=1,\ldots,E_h$ is a linear polynomial $u_{ih}\in\mathbb{P}_1(K_i)$
which can be expressed in terms of basis functions $\varphi_{ik}
\in\mathbb{P}_1(K_i)$ thus:
\beq\label{uhdef}
u_{ih}=\sum_{k=0}^{d}u_{ik}\varphi_{ik},\qquad i=1,\ldots,E_h.
\eeq
In this work, we use the modal {\it Taylor basis} which
 is defined by \cite{dglim,dglim2,luo2007}
\beq\label{taylor}
\varphi_{i0}\equiv1,\qquad \varphi_{ik}(\mathbf{x})=\frac{\mathbf{e}_k
  \cdot(\mathbf{x}-\mathbf{x}_{i0})}{\Delta x_{ik}},
\qquad k=1,\ldots, d,
\eeq
where $\mathbf{e}_k=(\delta_{kl})_{l=1}^d$ is a standard basis vector of
the Euclidean space $\mathbb{R}^d$,
$\mathbf{x}_{i0}=\frac{1}{|K_i|}\int_{K_i}\mathbf{x}\dx$ is the
centroid of element $K_i$ and $\Delta x_{ik}=\max_{\mathbf{x},\mathbf{y}\in K_i}
|\mathbf{e}_k\cdot(\mathbf{x}-\mathbf{y})|$ is a scaling factor
which improves the condition number of the mass matrix. The
Taylor degrees of freedom 
\beq
u_{i0}=\frac{1}{|K_i|}\int_{K_i}u_h(\mathbf{x})\dx,\qquad
u_{ik}=\Delta x_{ik}\pd{u_{ih}}{x_k}\Big|_{K_i}, \qquad k=1,\ldots, d
\eeq
represent the cell average and the constant partial derivatives of the
linear Taylor polynomial $u_{ih}$.

For simplicity, we assume that the whole boundary of $K_i$ lies in the
interior of the computational domain $\Omega_h=\bigcup_{i=1}^{E_h}K_i$ or
periodic boundary conditions are imposed. 
Let $\mathcal E_i$ denote the integer set containing the numbers of mesh cells
that share a side (boundary point in 1D, edge in 2D, face in
3D) $S_{ij}=\partial K_i\cap \partial K_j$ with $K_i$.
Substituting \eqref{uhdef} into a weak form of \eqref{ibvp-pde}
and using the notation $\mathbf{n}_{ij}$ for the unit outward
normal to
$S_{ij}$, we obtain $N_h=(d+1)E_h$ semi-discrete equations of the form
\beq\label{dg1high}
\sum_{l=0}^d\td{u_{il}}{t}\int_{K_i}\varphi_{ik}\varphi_{il}\dx=
-\sum_{j\in\mathcal E_i}\int_{S_{ij}}\varphi_{ik}
H(u_{ih},u_{jh},\mathbf{n}_{ij})\ds
+\int_{K_i}\nabla\varphi_{ik}\cdot\mathbf{f}(u_{ih})\dx,
\eeq
where $H(u_L,u_R,\mathbf{n})$
is a Lipschitz-continuous numerical approximation to the flux $\mathbf{f}(u)
\cdot\mathbf{n}$ such that  $H(u,u,\mathbf{n})=\mathbf{f}(u)\cdot\mathbf{n}$
and $H(u_L,u_R,\mathbf{n})+H(u_R,u_L,-\mathbf{n})=0$.

The system of equations \eqref{dg1high} for the time-dependent Taylor
degrees of freedom can be written as
\beq\label{dg1highw}
\sum_{l=0}^dm_{i,kl}\td{u_{il}}{t}=-\sum_{j\in\mathcal E_i}|S_{ij}|F_{ij,k}
+\int_{K_i}\nabla\varphi_{ik}\cdot\mathbf{f}(u_{ih})\dx,
\qquad i=1,\ldots,E_h,
\eeq
\beq\label{massTaylor}
m_{i,kl}=\int_{K_i}\varphi_{ik}\varphi_{il}\dx,
\qquad F_{ij,k}=\frac{1}{|S_{ij}|}\int_{S_{ij}}\varphi_{ik}
H(u_{ih},u_{jh},\mathbf{n}_{ij})\ds.
\eeq
Note that $m_{i,00}=|K_i|$ and $m_{i,0l}=0=m_{i,k0}$ for $k,l\in\{1,\ldots,d\}$.
If the Taylor basis is orthogonal, which is the case, e.g., on
uniform rectangular meshes, then the mass matrix $(m_{i,kl})_{k,l=0}^d$ is
diagonal.

For a given convex entropy $\eta(u)$, the corresponding entropy
variable $v=\eta'(u)$ can be approximated by the linear polynomial
$v_{ih}=\sum_{k=0}^dv_{ik}\varphi_{ik}$ with Taylor coefficients
\beq\label{vihdef}
v_{i0}=v(u_{i0}),\qquad v_{ik}=\eta''(u_{i0})u_{ik},
\qquad k=1,\ldots,d.
 \eeq
Let $b_i:\mathbb{P}_1(K_i)\times\mathbb{P}_1(K_i)
\to\mathbb{R}_0^+$ be a symmetric positive definite bilinear form. Define
\beq\label{ddef}
D_i(v_{ih},w_{ih})=b_i(v_{ih}-v_{i0},w_{ih}-w_{i0}),\qquad
v_{ih},w_{ih}\in \mathbb{P}_1(K_i).
\eeq
The assumption of positive definiteness (coercivity) implies that
$D_i(v_{i0},v_{i0})=0$ for the constant function $v_{ih}\equiv v_{i0}$
and $D_i(v_{ih},v_{ih})>0$ for all $v_{ih}\in\mathbb{P}_1(K_i)
\backslash\{v_{i0}\}$.

As noticed by Abgrall \cite{abgrall}, a dissipative
stabilization term of the form $\nu_iD_i(\varphi_{ik},v_{ih})$ can
be used to control the rate entropy production inside $K_i$.
Leaving the value of the stabilization parameter $\nu_i\ge 0$
unspecified for the time being, we consider the
entropy-corrected DG approximation
\beq\label{dg1ecorr}
\sum_{l=0}^dm_{i,kl}\td{u_{il}}{t}=
-\sum_{j\in\mathcal E_i}|S_{ij}|F_{ij,k}+
\int_{K_i}\nabla\varphi_{ik}\cdot\mathbf{f}(u_{ih})\dx
-\nu_iD_i(\varphi_{ik},v_{ih}).
\eeq
If $\varphi_{ik}$ are defined by \eqref{taylor},
then $D_i(\varphi_{i0},v_{ih})=0$ and $D_i(\varphi_{ik},v_{ih})=
b_i(\varphi_{ik},v_{ih}-v_{i0})$ for $k=1,\ldots,d$.

By the chain rule, we have $\pd{\eta(u_{ih})}{t}=\eta'(u_{ih})
\pd{u_{ih}}{t}$, where $\eta'(u_{ih})=v(u_{ih})$. It follows that
\beq\label{dnudt}
\td{}{t}\int_{K_i}\eta(u_{ih})\dx=
\int_{K_i}v(u_{ih})\pd{u_{ih}}{t}\dx
=\int_{K_i}\left[v_{ih}\pd{u_{ih}}{t}+
(v(u_{ih})-v_{ih})\pd{u_{ih}}{t}\right]\dx.
\eeq
Multiplying the semi-discrete conservation law
\eqref{dg1ecorr} by the Taylor coefficients $v_{ik}$
of $v_{ih}$, summing over $k=0,\ldots,d$,
and substituting the result into \eqref{dnudt},
we obtain the evolution equation
\beq
|K_i|\td{\eta_i}{t}=P_i(v_{ih},u_{ih})-\nu_iD_i(v_{ih},v_{ih}),
\eeq
where $\eta_i=\frac{1}{|K_i|}\int_{K_i}\eta(u_{ih}(\mathbf{x}))\dx$ is the
average entropy in $K_i$ and
\beq
P_i(v_{ih},u_{ih})=
-\sum_{j\in\mathcal E_i}|S_{ij}|\left(\sum_{k=0}^dv_{ik}F_{ij,k}\right)
+\int_{K_i}\left[\nabla v_{ih}\cdot\mathbf{f}(u_{ih})+
(v(u_{ih})-v_{ih})\pd{u_{ih}}{t}\right]\dx.
\eeq

\begin{rmk}
 The contribution of $(v(u_{ih})-v_{ih})\pd{u_{ih}}{t}$ vanishes for the
 square entropy $\eta(u)=\frac{u^2}{2}$. In general, the coefficients
 of $\pd{u_{ih}}{t}=\sum_{l=0}^d\td{u_{il}}{t}\varphi_{il}$ are obtained
 by solving system \eqref{dg1high} for $\td{u_{il}}{t}$.
\end{rmk}

Suppose that the numerical
fluxes $H(\cdot,\cdot,\cdot)$ and stabilization
 parameters $\nu_i\ge 0$
are chosen to satisfy 
\beq\label{eprodmax}
P_i(v_{ih},u_{ih})+\sum_{j\in\mathcal E_i}|S_{ij}|G_{ij}\le\nu_iD_i(v_{ih},v_{ih}),
\qquad i=1,\ldots,E_h,
\eeq
where $G_{ij}$ is a consistent approximation to the average value
$\frac{1}{|S_{ij}|}\int_{S_{ij}}\mathbf{q}(u)\cdot\mathbf{n}\ds$ of the entropy flux
$\mathbf{q}(u)$ associated with  $v(u)$. Then a
solution $u_{ih}$ of \eqref{dg1ecorr}
satisfies the cell entropy inequality
\beq\label{cellent}
|K_i|\td{\eta_i}{t}+\sum_{j\in\mathcal E_i}|S_{ij}|G_{ij}
\le 0
\eeq
which implies entropy stability of the semi-discrete
DG scheme. If $v_{ih}\ne v_{i0}$ then
$D_i(v_{ih},v_{ih})>0$ and for any finite value of
$P_i(v_{ih},u_{ih})$ the validity of
\eqref{eprodmax} can be readily enforced using
\beq\label{nudef}
\nu_i=\begin{cases}
\frac{\max\left\{0,P_i(v_{ih},u_{ih})+\sum_{j\in\mathcal E_i}|S_{ij}|
  G_{ij}\right\}}{D_i(v_{ih},v_{ih})} & \mbox{if}\ D_i(v_{ih},v_{ih})>0,\\
0 & \mbox{if}\ D_i(v_{ih},v_{ih})=0.
\end{cases}
\eeq
To ensure the validity of \eqref{eprodmax} for $v_{ih}=v_{i0}$ and
$D_i(v_{ih},v_{ih})=0$, we constrain the averaged fluxes
$F_{ij,0}=\frac{1}{|S_{ij}|}\int_{S_{ij}}
H(u_{ih},u_{jh},\mathbf{n}_{ij})\ds$ to
satisfy the entropy stability condition (cf. \cite{chen,entropyCG,pazner,tadmor})
\beq\label{condES}
(v_{j0}-v_{i0})F_{ij,0}\le
(\boldsymbol{\psi}(u_{j0})
-\boldsymbol{\psi}(u_{i0}))\cdot\mathbf{n}_{ij},
\eeq
where 
\beq
  \boldsymbol{\psi}(u)=v(u)\mathbf{f}(u)-\mathbf{q}(u).
  \eeq
  Condition  \eqref{condES} provides a useful tool for derivation of
  entropy stable
schemes. The following theorem shows that it is, indeed,
sufficient for \eqref{eprodmax} and \eqref{cellent} to hold
if $v_{ih}=v_{i0}$ and $D_i(v_{ih},v_{ih})=0$.

\begin{thm}[Entropy stability criterion]\label{thm1}
 The semi-discrete DG scheme \eqref{dg1ecorr} is
 entropy stable if the fluxes  $F_{ij,0}$ satisfy condition \eqref{condES}
 and the parameter $\nu_i$ is defined by \eqref{nudef}.
\end{thm} 

\begin{pf}
  For $D_i(v_{ih},v_{ih})>0$, the validity of \eqref{eprodmax} and
  \eqref{cellent} follows from \eqref{nudef}. If
  $D_i(v_{ih},v_{ih})=0$, then $v_{ih}=v_{i0}$
  and the entropy production term reduces to
  $P_i(v_{i0},u_{ih})=-\sum_{j\in\mathcal E_i}|S_{ij}|v_{i0}F_{ij,0}$.
  In view of the assumption that the flux $F_{ij,0}$ satisfies
  condition  \eqref{condES}, we have the estimate
$$
  v_{i0}F_{ij,0}=\frac{v_{j0}+v_{i0}}{2}F_{ij,0}-
  \frac{v_{j0}-v_{i0}}{2}F_{ij,0}\ge \frac{v_{j0}+v_{i0}}{2}F_{ij,0}-
  \frac12(\boldsymbol{\psi}(u_{j0})
-\boldsymbol{\psi}(u_{i0}))\cdot\mathbf{n}_{ij}.
$$
Exploiting the fact that
$\boldsymbol{\psi}(u_{i0})\cdot\sum_{j\in\mathcal E_i}|S_{ij}|\mathbf{n}_{ij}=0$,
we can now estimate $P_i$ as follows:
\begin{align*}
  P_i(v_{i0},u_{ih})&\le -\sum_{j\in\mathcal E_i}|S_{ij}|
  \left(\frac{v_{j0}+v_{i0}}{2}F_{ij,0}-
  \frac12(\boldsymbol{\psi}(u_{j0})
-\boldsymbol{\psi}(u_{i0}))\cdot\mathbf{n}_{ij}
\right)\\
&=-\sum_{j\in\mathcal E_i}|S_{ij}|
  \left(\frac{v_{j0}+v_{i0}}{2}F_{ij,0}-
  \frac12(\boldsymbol{\psi}(u_{j0})
+\boldsymbol{\psi}(u_{i0}))\cdot\mathbf{n}_{ij}
\right).
  \end{align*}
Thus, we have $P_i(v_{i0},u_{ih})+\sum_{j\in\mathcal E_i}|S_{ij}|G_{ij}\le 0
=\nu_iD(v_{i0},v_{i0})$ for $G_{ij}$ defined by
$$
G_{ij}=\frac{v_{j0}+v_{i0}}{2}F_{ij,0}-
  \frac12(\boldsymbol{\psi}(u_{j0})
+\boldsymbol{\psi}(u_{i0}))\cdot\mathbf{n}_{ij}.
$$
This proves the validity of inequalities \eqref{eprodmax} and
  \eqref{cellent} in the case $v_{ih}=v_{i0}$.
   \qquad\proofend
\end{pf}

\section{Flux limiting and entropy corrections}
\label{sec:FL}

The above analysis provides general guidelines for the design
of entropy stable DG-$\mathbb{P}_1$ schemes. In this section, we present
practical algorithms for calculating numerical fluxes
that ensure not only entropy stability but also positivity
preservation for the cell averages. We also define a stabilization
operator $D_i(\cdot,v_{ih})$ which penalizes the gradients of $u_{ih}$
to satisfy condition \eqref{eprodmax}. The methodology to be presented combines
and generalizes
the entropy correction tools developed in \cite{abgrall,ranocha}
and \cite{entropyCG}.

Many nonlinear high-resolution schemes are based on the idea of
blending a property-preserving low-order flux $F_{ij}^L$ and a
high-order {\it target flux} $F_{ij}^H$. The former is supposed
to satisfy inequality constraints that guarantee entropy
stability and/or the validity of (local) maximum principles.
Using an adaptively chosen correction factor $\alpha_{ij}\in[0,1]$,
the convex combination $F_{ij}=(1-\alpha_{ij})F_{ij}^L+\alpha_{ij}F_{ij}^H$
can be constrained to satisfy them as well. If the flux $F_{ij}^H$
possesses the desired properties, then $\alpha_{ij}=1$ is the
optimal choice. Otherwise, the best inequality-constrained
approximation corresponds to the largest value of $\alpha_{ij}\in[0,1)$
for which the property constraints can be shown to hold. This design
philosophy is common to flux-corrected transport (FCT) algorithms,
total variation diminishing (TVD) methods, and other representatives
of local extremum diminishing (LED) flux correction schemes
\cite{afc1}. The possibility of using flux limiting to enforce entropy
inequalities in addition to LED constraints  was explored
in \cite{entropyCG} in the context of continuous Galerkin methods
and Lagrange finite elements. Building on this work, we equip
our Taylor DG-$\mathbb{P}_1$ scheme \eqref{dg1ecorr} with limited
fluxes of the form
\beq\label{fluxlim}
F_{ij,k}=(1-\alpha_{ij})F_{ij,k}^L+\alpha_{ij}F_{ij,k}^H,\qquad
k=0,\ldots,d.
\eeq
The general flux function of the local Lax-Friedrichs (LLF) method is
defined by (cf. \cite{chen,Guermond2016})
 \beq 
  H_{\rm LLF}(u_L,u_R,\mathbf{n})=
  \frac{\mathbf{f}(u_R)+\mathbf{f}(u_L)}{2}\cdot\mathbf{n}
  -\lambda(u_L,u_R,\mathbf{n})
  \frac{u_R-u_L}{2},
  \eeq
  where   $\lambda(u_L,u_R,\mathbf{n})$ is the maximal speed of wave 
  propagation in the direction parallel to $\mathbf{n}$, i.e.,
  \beq
  \lambda(u_L,u_R,\mathbf{n})=\max_{\omega\in[0,1]}|
\mathbf{f}'(\omega u_L+(1-\omega) u_R)\cdot\mathbf{n}|.
\eeq
  
As shown by Chen and Shu \cite{chen}, 
condition \eqref{condES} is satisfied for the low-order
LLF flux
\beq
F_{ij,k}^L=\frac{1}{|S_{ij}|}\int_{S_{ij}}\varphi_{ik}
H_{\rm LLF}(u_{i0},u_{j0},\mathbf{n}_{ij})\dx,\qquad k=0,\ldots,d.
\eeq
The high-order LLF flux $F^H_{ij,k}$ of the unconstrained
DG-$\mathbb{P}_1$ approximation is given by
\beq
F_{ij,k}^H=\frac{1}{|S_{ij}|}\int_{S_{ij}}\varphi_{ik}
H_{\rm LLF}(u_{ih},u_{jh},\mathbf{n}_{ij})\dx,\qquad k=0,\ldots,d.
\eeq

Let us first determine an entropy correction factor $\alpha_{ij}^{\rm ES}$
such that the limited LLF flux
\beq\label{fluxlim0}
F_{ij,0}=(1-\alpha_{ij})F_{ij,0}^L+\alpha_{ij}F_{ij,0}^H
\eeq
satisfies \eqref{condES} for any $\alpha_{ij}\in[0,\alpha_{ij}^{\rm ES}]$. 
Following the derivation of algebraic entropy fixes for continuous
Galerkin methods in
\cite{entropyCG}, we use $\alpha_{ij}^{\rm ES}\in\{\alpha_{ij}^{\rm ES1},
\alpha_{ij}^{\rm ES2},\alpha_{ij}^{\rm ES3}\}$, where $\alpha_{ij}^{\rm ES1}
\le\alpha_{ij}^{\rm ES2}\le\alpha_{ij}^{\rm ES3}$ are correction factors
corresponding to different definitions of the  nonnegative
bound $Q_{ij}$ in the 
formula
\beq\label{alphaES}
\alpha_{ij}^{\rm ES}=
\begin{cases}
  \frac{Q_{ij}}{P_{ij}} & \mbox{if}\ P_{ij}>Q_{ij},\\
 1 & \mbox{otherwise},
\end{cases}  \qquad
P_{ij}=(v_{j0}-v_{i0})(F_{ij,0}^H-F_{ij,0}^L).
\eeq
The least dissipative entropy stability preserving
upper bound for the rate $\alpha_{ij}P_{ij}$ of entropy increase
in $K_i$ due to the limited antidiffusive flux
$\alpha_{ij}(F_{ij,0}^H-F_{ij,0}^L)$
is given by
\begin{align}
Q_{ij}^{\rm ES1}&=(\boldsymbol{\psi}(u_{j0})
  -\boldsymbol{\psi}(u_{i0}))\cdot\mathbf{n}_{ij}
  -(v_{j0}-v_{i0})F_{ij,0}^L =(v_{j0}-v_{i0})
 \frac{\lambda_{ij}}{2}(u_{j0}-u_{i0})\nonumber\\
&+(\boldsymbol{\psi}(u_{j0})
-\boldsymbol{\psi}(u_{i0}))\cdot\mathbf{n}_{ij}-(v_{j0}-v_{i0})\frac{\mathbf{f}(u_{j0})+\mathbf{f}(u_{i0})}{2}\cdot\mathbf{n}_{ij}.
  \end{align}
  Nonnegativity of $Q_{ij}^{\rm ES1}$ follows from the fact that 
  condition \eqref{condES} holds for the low-order flux $F_{ij,0}^L$.
  It is easy to verify that $\alpha_{ij}P_{ij}\le 
  Q_{ij}^{\rm ES1}$ is a
  sufficient condition for \eqref{fluxlim0} to satisfy \eqref{condES}.
  In accordance with Tadmor's \cite{tadmor,tadmor2016} comparison principle
  for entropy conservative and entropy stable schemes, definition 
  \eqref{alphaES} guarantees entropy stability for any $Q_{ij}\in[0,
  Q_{ij}^{\rm ES1}]$. To prevent the limited
  flux \eqref{fluxlim0} from producing more entropy than the
   centered flux
  $F_{ij,0}^{\rm CD}=\frac{\mathbf{f}(u_{j0})
    +\mathbf{f}(u_{i0})}{2}\cdot\mathbf{n}_{ij}$, we use
  \beq
  Q_{ij}^{\rm ES2}=\max\left\{0,(v_{j0}-v_{i0})
 \frac{\lambda_{ij}}{2}(u_{j0}-u_{i0})+
  \min\left\{0,Q_{ij}^{\rm CD}\right\}\right\},
  \eeq
  where $\lambda_{ij}=\lambda(u_{i0},u_{j0},\mathbf{n}_{ij})$
 and $Q_{ij}^{\rm CD}=(\boldsymbol{\psi}(u_{j0})
  -\boldsymbol{\psi}(u_{i0}))\cdot\mathbf{n}_{ij}
  -(v_{j0}-v_{i0})F_{ij,0}^{\rm CD}$. Definition
   \beq
   Q_{ij}^{\rm ES3}=\max\left\{0,(v_{j0}-v_{i0})\left(
 \frac{\lambda_{ij}}{2}-\nu_{ij}\right)(u_{j0}-u_{i0}) +
  \min\left\{0,Q_{ij}^{\rm CD}\right\}\right\}
  \eeq
  makes it possible to further increase the levels of
   entropy dissipation using (cf. \cite{entropyCG})
  \beq
  \nu_{ij}=\begin{cases}
  \frac{
    \max\left\{0,\left[\frac{\mathbf{f}(u_{j0})+\mathbf{f}(u_{i0})}{2}
    -\mathbf{f}\left(\frac{u_{j0}+u_{i0}}{2}\right)\right]
    \cdot\mathbf{n}_{ij}
    \right\}}{v_{j0}-v_{i0}} & \mbox{if}\ v_{j0}\ne v_{i0},\\
 0 & \mbox{if}\ v_{j0}\ne v_{i0}.
\end{cases}
\eeq
The numerical examples of Section \ref{sec:num} illustrate the
ramifications of different choices of $Q_{ij}$.

It is also possible to find a correction factor $\alpha_{ij}^{\rm BP}$
which ensures preservation of the local bounds 
\beq\label{lbounds}
u_{i0}^{\min}=\min_{j\in\mathcal E_i^*}u_{j0},\qquad
u_{i0}^{\max}:=\max_{j\in\mathcal E_i^*}u_{j0},
\eeq
where $\mathcal E_i^*$ is the set containing the numbers of all elements
$E_j$ such that $E_i$ and $E_j$ have at least one common vertex. Note that
$i=j$ is also an element of $\mathcal E_i^*$.
Adapting the {\it monolithic
  convex limiting} (MCL) strategy \cite{convex,EG-MCL,entropyCG} to the
DG setting, we impose the inequality constraints
\beq\label{afc-bounds}
u_{i0}^{\min}
\le \bar u_{ij,0}^*:=
\bar u_{ij,0}+\frac{\alpha_{ij}F_{ij,0}^A}{\lambda_{ij}}
\le u_{i0}^{\max},
\eeq
where $F_{ij,0}^A=F_{ij,0}^H-F_{ij,0}^L$ is the raw antidiffusive flux and
$\bar u_{ij,0}$ are intermediate states such that
\beq
\min\{u_{i0},u_{j0}\}\le
\bar u_{ij,0}:=\frac{u_{j0}+u_{i0}}{2}-\frac{\mathbf{f}(u_{j0})-\mathbf{f}(u_{j0})
}{2\lambda_{ij}}\cdot\mathbf{n}_{ij}\le\max\{u_{i0},u_{j0}\}.
\eeq
The validity of the MCL constraints \eqref{afc-bounds}
is guaranteed for $\alpha_{ij}
\le\alpha_{ij}^{\rm BP}$, where (cf. \cite{convex,EG-MCL})
\beq\label{alphaBP}
\alpha_{ij}^{\rm BP}=\begin{cases}
\min\,\left\{1,
\frac{\lambda_{ij}\min\,
  \{u_{i0}^{\max}-\bar u_{ij,0},\bar u_{ji,0}-u_{j0}^{\min}\}}{F_{ij,0}^A}
\right\}
  & \mbox{if}\  F_{ij,0}^A>0,
\\[0.25cm]
 1 & \mbox{if}\  F_{ij,0}^A=0,
\\[0.25cm]
\min\,\left\{1,
\frac{\lambda_{ij}\min\,
  \{u_{i0}^{\min}-\bar u_{ij,0},\bar u_{ji,0}-u_{j0}^{\max}\}}{F_{ij,0}^A}
\right\}
  & \mbox{if}\  F_{ij,0}^A<0.
\end{cases}
\eeq
In the next section, we prove that the fully discrete version of
the flux-limited DG scheme satisfies a local maximum principle
if integration in time is performed using an SSP Runge-Kutta
method.

In light of the above, the limited fluxes \eqref{fluxlim} should be defined
using $\alpha_{ij}=\alpha_{ij}^{\rm ES}$ to enforce entropy stability,
$\alpha_{ij}=\alpha_{ij}^{\rm BP}$ to ensure preservation of local
bounds, and $\alpha_{ij}=\min\{\alpha_{ij}^{\rm ES},\alpha_{ij}^{\rm BP}\}$
to provide both properties.
To conclude the derivation of the semi-discrete property-preserving
DG scheme, we need to choose a coercive bilinear form $b_i(\cdot,\cdot)$
for definition \eqref{ddef} of the entropy stabilization operator $D_i$.
The bilinear form of the Rusanov-type dissipation term employed in
\cite{ranocha} can be written as
\beq
b_i(v,w)=\int_{K_i}vw\dx,\qquad v,w\in L^2(K_i)
\eeq
and induces
\beq\label{DiRus}
D_i(\varphi_{ik},v_{ih})=\int_{K_i}\varphi_{ik}(v_{ih}-v_{i0})\dx,
\qquad k=0,\ldots,d,
\eeq
where $v_{i0}$ is the cell average of $v_{ih}$. This representation of
$D_i(\varphi_{ik},v_{ih})$ is suitable not only for nodal (e.g., Lagrange,
Bernstein or Gauss-Lobatto) finite element bases but also for the Taylor
basis \eqref{taylor}.

In essence, the addition of $\nu_iD_i(\varphi_{ik},v_{ih})$ on the right-hand
side of \eqref{dg1ecorr} penalizes the deviations of
$u_{ih}(\mathbf{x})=u_{i0}+(\mathbf{x}-
\mathbf{x}_{i0})\cdot\nabla u_{ih}$ from $u_{i0}=u_{ih}(\mathbf{x}_{i0})$.
For the Taylor basis \eqref{taylor}, we have
\beq\label{DiTaylor}
D_i(\varphi_{i0},v_{ih})=0,\qquad D_i(\varphi_{ik},v_{ih})
=\sum_{l=1}^dm_{i,kl}v_{il},\qquad k=1,\ldots,d,
\eeq
where $m_{i,kl}$ are entries of the consistent Taylor mass matrix, 
as defined by \eqref{massTaylor}. For an orthogonal Taylor basis,
the mass matrix is diagonal and, therefore, definition \eqref{DiTaylor}
simplifies to
\beq\label{DiSimp}
D_i(\varphi_{ik},v_{ih})=m_{i,kk}v_{ik},\quad k=1,\ldots,d.
\eeq
This simplified definition may be used in the case of a non-orthogonal
Taylor basis as well. The associated coercive bilinear form is given
by $\tilde b_{ih}(v_{ih},w_{ih})=\sum_{k=0}^dm_{i,kk}v_{ik}w_{ik}$.

\begin{rmk}
  For a nodal basis $\{\hat\varphi_{ik}\}_{k=1}^{d+1}$, the
  $L^2$ scalar product $b_i(v_{ih},w_{ih})$ 
  can be approximated by $\tilde b_{ih}(v_{ih},w_{ih})=\sum_{k=1}^{d+1}\hat
  m_{i,k}\hat v_{ik}\hat w_{ik}$, where $\hat v_{ik}$ are nodal values
  or Bernstein coefficients and $\hat m_{i,k}$ are diagonal entries of the
lumped element mass matrix. This approach corresponds to inexact
numerical integration for $b_i(v_{ih},w_{ih})$ and replaces \eqref{DiRus}
with $D_i(\hat\varphi_{ik},v_{ih})=\hat m_{i,k}(\hat v_{ik}-v_{i0})$,
cf. \cite{Abgrall2006,abgrall}.
\end{rmk}

\begin{rmk}
  Algebraic counterparts of the Rusanov dissipation operator
   \eqref{DiRus} are
  often used to construct low-order components of
  nonlinear high-resolution finite element schemes
  \cite{Abgrall2006,Guermond2014,Guermond2016,RD-BFCT,afclps}.
  Different authors write them in different forms and give them
  different names. Representation of $D_i(\cdot,\cdot)$
  in the basis-independent form \eqref{DiRus} was used in \cite{afclps}
  in the context of algebraic flux correction (AFC) schemes for continuous
  finite elements. We refer the interested reader to \cite{afclps} for an
  in-depth discussion of its properties and existing
  relationships to other forms of Rusanov dissipation.
 \end{rmk}

\section{Time discretization and slope limiting}
\label{sec:SL}

Let $F_{ij,k}$ be the limited fluxes defined by \eqref{fluxlim} and
$\nu_i$ the stabilization parameter defined by \eqref{nudef}.
Substituting the diagonal form \eqref{DiSimp} of $D_i(\varphi_{ik},v_{ih})$
into \eqref{dg1ecorr}, we obtain the semi-discrete problem
\beq\label{dg1ecorr2}
m_{i,kk}\td{u_{ik}}{t}=
-\sum_{j\in\mathcal E_i}|S_{ij}|F_{ij,k}+
\int_{K_i}\nabla\varphi_{ik}\cdot\mathbf{f}(u_{ih})\dx
-(1-\delta_{0k})\left[\nu_im_{i,kk}v_{ik}+\sum_{l=1\atop l\ne k}^dm_{i,kl}\dot u_{i,l}\right],
\eeq
where $\dot u_{i,l}$ is the time derivative of the Taylor degree of
freedom $u_{i,l}$, as defined by \eqref{dg1ecorr}, or a suitable
approximation thereof. For simplicity, we set $\dot u_{i,l}:=0$
and consider the reduced system
\begin{align} \label{dg1ecorr3}
|K_i|\td{u_{i0}}{t}
&=-\sum_{j\in\mathcal E_i}|S_{ij}|F_{ij,0},\\
m_{i,kk}\left(\td{u_{ik}}{t}+\nu_iv_{ik}\right)
&=-\sum_{j\in\mathcal E_i}|S_{ij}|F_{ij,k}+
\int_{K_i}\nabla\varphi_{ik}\cdot\mathbf{f}(u_{ih})\dx,\qquad k=1,\ldots,d
 \label{dg1ecorr4}
\end{align}
which is equivalent to \eqref{dg1ecorr2} in the case of an orthogonal
Taylor basis and corresponds to a locally conservative lumped-mass
approximation otherwise (see \cite{dglim2} for a discussion of
mass lumping and appropriate treatment of time derivatives in Taylor basis
DG methods).

We advance $u_{i0}$ in time using an explicit SSP
Runge-Kutta method with at least two stages (for accuracy and
stability reasons). Each stage has the structure of the
forward Euler update
\beq\label{FEupdate}
|K_i|\bar u_{i0}=|K_i|u_{i0}-\Delta t\sum_{j\in\mathcal E_i}|S_{ij}|F_{ij,0},
\eeq
where $\Delta t>0$ is the time step and the numerical flux $F_{ij,0}$ is
defined at the same old time level as $u_{i0}$.

The zero-order term $m_{i,kk}\nu_iv_{ik}$ on the left-hand side of
equation \eqref{dg1ecorr4} must be treated implicitly to
avoid uncontrolled growth of the solution gradients (in the absence
of limiting) and severe time step restrictions. Approximating
$v_{ik}=\eta''(u_{i0})u_{ik},\ k=1,\ldots,d$ by $\eta''(u_{i0})\bar u_{ik}$,
we update $u_{ik}$ as follows:
\beq\label{FEupdate2}
m_{i,kk}\left(1+\Delta t\nu_i\eta''(u_{i0})\right)\bar u_{ik}
=m_{i,kk}u_{ik}-\Delta t\left(
\sum_{j\in\mathcal E_i}|S_{ij}|F_{ij,k}-
\int_{K_i}\nabla\varphi_{ik}\cdot\mathbf{f}(u_{ih})\dx\right).
\eeq
Note that $\eta''(u_{i0})\ge 0$ for any convex entropy $\eta(u)$. 
For the square entropy $\eta(u)=\frac{u^2}{2}$, we have $\eta''(u_{i0})=1$, and
the coefficient of the left-hand side reduces to $m_{i,kk}+\Delta t\nu_i$. 
All quantities that appear on the right-hand side, as well as the value
of $\nu_i\ge 0$ are defined using the old solution $u_h$.

The semi-implicit treatment of $\nu_iv_{ik}$ corresponds to Patankar-type
source term linearization \cite{burchard,patankar} for positivity-preserving 
schemes. It guarantees that the magnitude
of $\nabla \bar u_{ih}$ decreases as the value of the entropy viscosity
coefficient $\nu_i$ is increased. Hence, the addition of this
term penalizes steep gradients in the same manner as it does
at the semi-discrete level. In our experience, the explicit 
treatment of $\nu_iv_{ik}$ has a devastating effect on
the entropy stability behavior of the fully discrete scheme.
Indeed, it may fail to control the magnitude of $\nabla \bar u_{ih}$
in the desired manner and is likely to increase it significantly
in situations when the rate of entropy production and the value
 of $\nu_i>0$ become very large.

 Importantly, the use of $\alpha_{ij}\le\alpha_{ij}^{\rm BP}$ in the formula
 for $F_{ij,0}$ guarantees preservation of invariant domains and local
 bounds for the cell averages of the entropy-stabilized DG-$\mathbb{P}_0$
 approximation. The following adapted version of a theorem proved in
 \cite{convex} is presented here for the reader's convenience.
 
\begin{thm}[Preservation of local bounds \cite{convex}]\label{thm2}
  Let $u_{i0}$ be evolved using \eqref{FEupdate}, where $F_{ij,0}$ is
  defined by \eqref{fluxlim0} with $\alpha_{ij}\in[0,\alpha_{ij}^{\rm BP}]$
  and $\alpha_{ij}^{\rm BP}$ defined by \eqref{alphaBP}.
  Choose the time step $\Delta t$ to satisfy
  \beq\label{cfl}
\Delta t\sum_{j\in\mathcal E_i}|S_{ij}|\lambda_{ij}\le |K_i|.
\eeq
Then the cell average $\bar u_{i0}$ satisfies the local maximum principle
\beq\label{lmp}
u_{i0}^{\min}\le\bar u_{i0}\le u_{i0}^{\max},
\eeq
where $u_{i0}^{\min}$ and $u_{i0}^{\max}$ are the local bounds defined in
\eqref{lbounds}.

\end{thm} 

\begin{pf}
  The SSP Runge-Kutta stage \eqref{FEupdate} can be written in the bar state
  form (cf. \cite{Guermond2016,convex})
  \beq
  \bar u_{i0}=u_{i0}+\frac{\Delta t}{|K_i|}\sum_{j\in\mathcal E_i}
  |S_{ij}|\lambda_{ij}(\bar u_{ij,0}^{*}-u_{i0})
=\left(1-\frac{\Delta t}{|K_i|}\sum_{j\in\mathcal E_i}
  |S_{ij}|\lambda_{ij}\right)u_{i0}+\frac{\Delta t}{|K_i|}\sum_{j\in\mathcal E_i}
  |S_{ij}|\lambda_{ij}\bar u_{ij,0}^{*},
  \eeq
  where $\bar u_{ij,0}^*$ are the bound-preserving flux-corrected
    bar states defined by \eqref{afc-bounds}. For time steps
    satisfying condition \eqref{cfl}, the result $\bar u_{i0}$ is a convex
    combination of  $u_{i0}\in[u_{i0}^{\min},u_{i0}^{\max}]$  and
 $\bar u_{ij,0}^*\in[u_{i0}^{\min},u_{i0}^{\max}]$
    , which
        proves the validity of estimate \eqref{lmp}.
   \qquad\proofend
\end{pf}
  A local maximum principle can also be proved for the
  steady-state limit of \eqref{FEupdate}, see Appendix of
  \cite{convex}.
  \smallskip
  
In addition to flux limiting, the slopes of the DG solution
can be adjusted to ensure that the value of $\bar u_{ih}$ at any point
$\mathbf{x}\in K_i$
will be bounded by the local minimum $\bar u_{i,0}^{\min}$ and maximum
$\bar u_{i,0}^{\max}$ of the property-preserving cell averages. Since
the linear polynomial $\bar u_{ih}$, whose Taylor coefficients
$\bar u_{ik}$ are defined by \eqref{FEupdate} and \eqref{FEupdate2},
attains its maxima and minima
at the vertices $\mathbf{x}_{i1},\ldots,\mathbf{x}_{iN}$
of $K_i$, we constrain its gradients using a
vertex-based version \cite{dglim,dglim2} of the Barth-Jespersen
slope limiter \cite{barthjesp}. Let $\mathcal E_{ip},\ p\in\{1,\ldots,N\}$
be a subset of
$\mathcal E_i^*$ containing the numbers of 
cells to which the vertex $\mathbf{x}_{p}$
belongs (including $K_i$). Multiplying $\bar u_{ik},
\ k=1,\ldots,d$ by a slope limiting factor $\beta_i\in[0,1]$, we constrain
\beq\label{slim}
\bar u_{ih}^*(\mathbf{x})=\bar u_{i0}+\beta_i\sum_{k=1}^d\bar u_{ik}
\varphi_{ik}(\mathbf{x})=\bar u_{i0}+\beta_i
(\mathbf{x}-\mathbf{x}_{i0})\cdot\nabla\bar u_{ih},\qquad \mathbf{x}\in K_i
\eeq
to satisfy the inequality constraints
\beq\label{vblmp}
\min_{j\in\mathcal E_{ip}}\bar u_{j0}=:\bar u_{ip}^{\min}
\le \bar u_{ih}^*(\mathbf{x}_{ip})
\le \bar u_{ip}^{\max}:=\max_{j\in\mathcal E_{ip}}\bar u_{j0},\qquad p=1,\ldots,N.
\eeq
The vertex-based slope limiter employed in \cite{dglim}
accomplishes this task by using the correction factor
\beq\label{betadef}
\beta_i=
\min_{1\le p \le N}\begin{cases}
  \min\left\{1,\frac{\bar u_{ip}^{\max}-\bar u_{i0}}{(\mathbf{x}_{ip}-\mathbf{x}_{i0})\cdot\nabla\bar u_{ih}}
    \right\} & \mbox{if}\  (\mathbf{x}_{ip}-\mathbf{x}_{i0})\cdot\nabla\bar u_{ih} >0,\\
           1 &\mbox{if}\ (\mathbf{x}_{ip}-\mathbf{x}_{i0})\cdot\nabla\bar u_{ih}=0,\\
           \min\left\{1,\frac{\bar u_{ip}^{\min}-\bar u_{i0}}{(\mathbf{x}_{ip}-\mathbf{x}_{i0})\cdot\nabla\bar u_{ih}}
    \right\} & \mbox{if}\  (\mathbf{x}_{ip}-\mathbf{x}_{i0})\cdot\nabla\bar u_{ih} <0.
    \end{cases}
\eeq
\begin{rmk}
  In most cases, the slope-limited version of a DG-$\mathbb{P}_1$
scheme  produces nonoscillatory
results even if the bound-preserving flux limiter is deactivated.
However, the assertion of Theorem \ref{thm2} is generally not true
for the high-order LLF fluxes $F_{ij,0}^H$ or the entropy-limited
fluxes $F_{ij,0}$ defined by \eqref{fluxlim0} with
$\alpha_{ij}=\alpha_{ij}^{\rm ES}$. As
noticed by Moe et al. \cite{Moe2017}, it is essential to use a
flux limiter for cell averages in addition to
slope limiting if positivity preservation is a must. Indeed,
the cell averages are the main unknowns of the discrete problem.
If they are constrained properly, slope limiting can be replaced
with less aggressive accuracy-preserving gradient corrections.
For example, a suitable  smoothness indicator (cf. \cite{kriv2004,persson})
can be used to
increase the value of the gradient penalization coefficient
$\nu_i$ in troubled cells. Following the analysis of
monolithic
algebraic flux correction schemes \cite{afc_analysis,CL-diss}, the
formula for $\nu_i$ can be designed to provide Lipschitz
continuity of $\nu_iD_i$, which is an essential
requirement for proving well-posedness
of the nonlinear discrete problem and convergence to
steady-state solutions. If necessary, the
final output may be postprocessed
using the vertex-based slope limiter to eliminate untershoots
and overshoots (if any) in DG-$\mathbb{P}_1$ solutions to be visualized or used
to calculate derived quantities.
\end{rmk}

\begin{rmk}
  The assumptions of Theorem \ref{thm1} do not guarantee that the cell
  averages defined by \eqref{FEupdate} will satisfy a fully discrete
  entropy inequality. However, an explicit SSP Runge-Kutta time discretization
  of the low-order LLF scheme corresponding
  to \eqref{fluxlim0} with $\alpha_{ij}=0$  is entropy stable for any entropy
  pair, as shown by Guermond and Popov \cite{Guermond2016} in the
  context of first-order continuous finite element approximations that
  exhibit the same structure \cite{Guermond2019}.
  Hence, the algorithm for calculating
  $\alpha_{ij}$ can be modified to check and enforce inequality
  constraints that imply fully discrete entropy stability. 
  \end{rmk}

\section{Numerical examples}\label{sec:num}

In this section, we perform numerical experiments for two nonlinear scalar
test problems in 2D. In the description of numerical results, the methods under
investigation are labeled as follows:
\begin{itemize}
\item DG0: standard DG-$\mathbb{P}_0$ discretization using the
  low-order LLF
  fluxes $F_{ij,k}^L$,
\item DG1: standard DG-$\mathbb{P}_1$ discretization using the
  high-order LLF
  fluxes $F_{ij,k}^H$,
\item ESX: entropy stable DG scheme using \eqref{fluxlim}
  with $\alpha_{ij}=\alpha_{ij}^{\rm ESX},\ X\in\{1,2,3\}$.
\end{itemize}  
We append the letters F and S to abbreviations of methods that use bound-preserving flux and slope limiting, respectively. For example, ES1F is the DG scheme
defined by \eqref{fluxlim} with $\alpha_{ij}=\min\{\alpha_{ij}^{\rm ES1},
\alpha_{ij}^{\rm BP}\}$ and ES1FS is its slope-limited counterpart which
adjusts the gradients using \eqref{slim} and \eqref{betadef}.

In all numerical examples, we use the square entropy $\eta(u)=\frac{u^2}2$
and the associated entropy variable $v(u)=u$. Computations are
performed on uniform rectangular meshes, on which the Taylor
basis is orthogonal. Numerical solutions are advanced in time using
the explicit third-order three-stage SSP Runge-Kutta method \cite{ssprev}
and time steps satisfying the CFL-like condition \eqref{cfl}
of Theorem \ref{thm2}. For visualization purposes, we project
DG solutions into the space of continuous bilinear elements.

\subsection{KPP problem}
The KPP problem \cite{Guermond2016,Guermond2017,kpp} is a
challenging nonlinear test for verification of entropy stability
properties. We use this problem to 
test different components of the method that we propose. 
In this series of 2D experiments,
we solve equation \eqref{ibvp-pde} with the nonlinear and
nonconvex flux function
\beq
\mathbf{f}(u)=(\sin(u),\cos(u))
\eeq
 in the computational domain
$\Omega_h=(-2,2)\times(-2.5,1.5)$ using the initial
condition
\beq
u_0(x,y)=\begin{cases}
\frac{14\pi}{4} & \mbox{if}\quad \sqrt{x^2+y^2}\le 1,\\
\frac{\pi}{4} & \mbox{otherwise}.
\end{cases}
\eeq
The entropy flux corresponding to $\eta(u)=\frac{u^2}2$ is
$\mathbf{q}(u)=(u\sin(u)+\cos(u),u\cos(u)-\sin(u))$.
A simple (but rather pessimistic) upper bound for the guaranteed
maximum speed (GMS) is $\lambda=1$. More accurate GMS estimates
can be found in \cite{Guermond2017}. The exact solution exhibits
a two-dimensional rotating wave structure, which is difficult to
capture in numerical simulations using high-order
methods. The main challenge of this test is to
prevent possible convergence to wrong weak solutions.

\begin{figure}
  
\begin{minipage}[t]{0.33\textwidth}
  \centering\small DG0, $u_h\in[0.785,10.992]$

\includegraphics[width=0.9\textwidth,trim=100 20 100 20,clip]{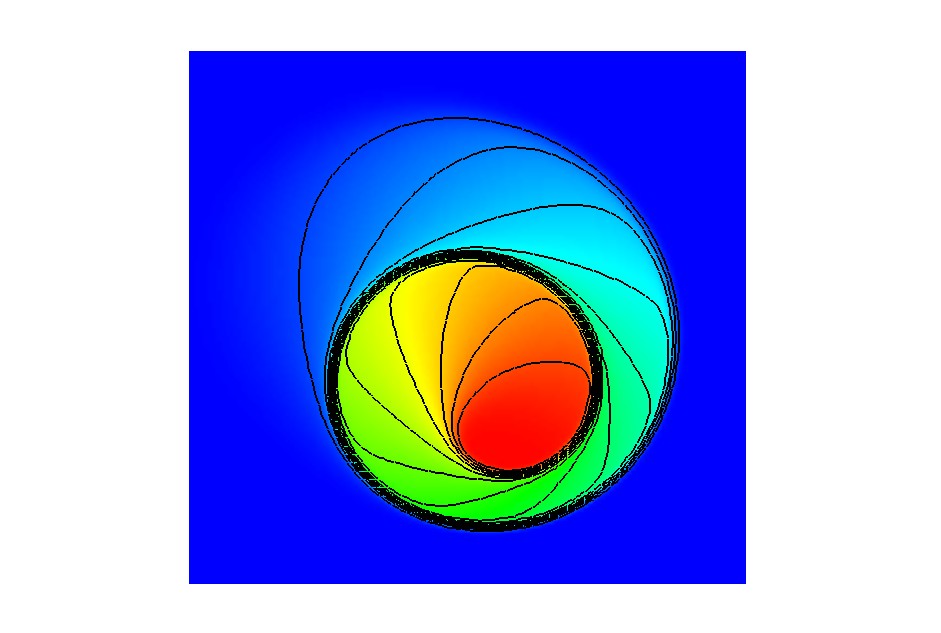}

\end{minipage}%
\begin{minipage}[t]{0.33\textwidth}

\centering\small  DG1, $u_h\in[-1.794,14.741]$

\includegraphics[width=0.9\textwidth,trim=100 20 100 20,clip]{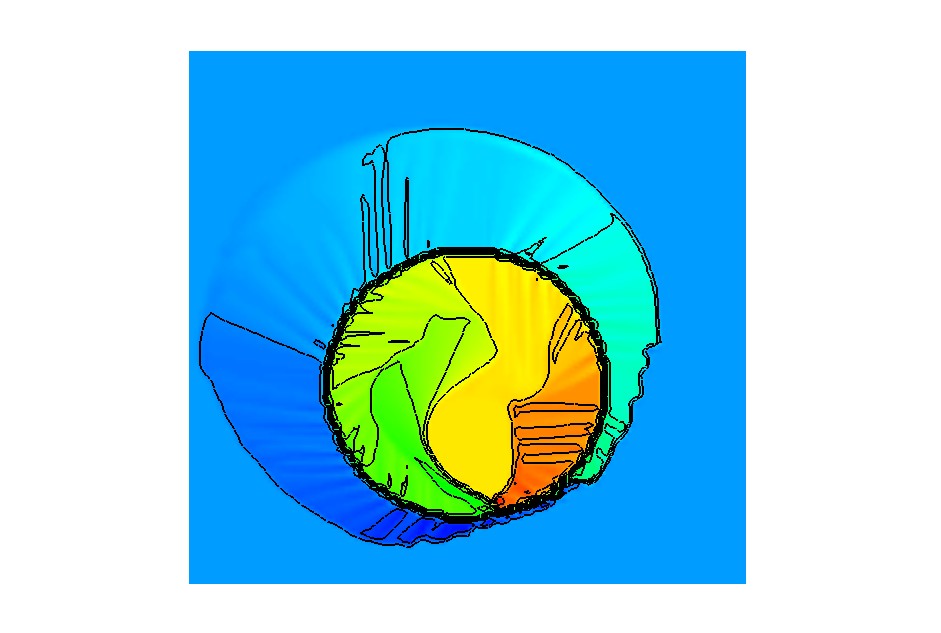}

\end{minipage}%
\begin{minipage}[t]{0.33\textwidth}

  \centering\small  DG1S, $u_h\in[0.785,10.996]$

\includegraphics[width=0.9\textwidth,trim=100 20 100 20,clip]{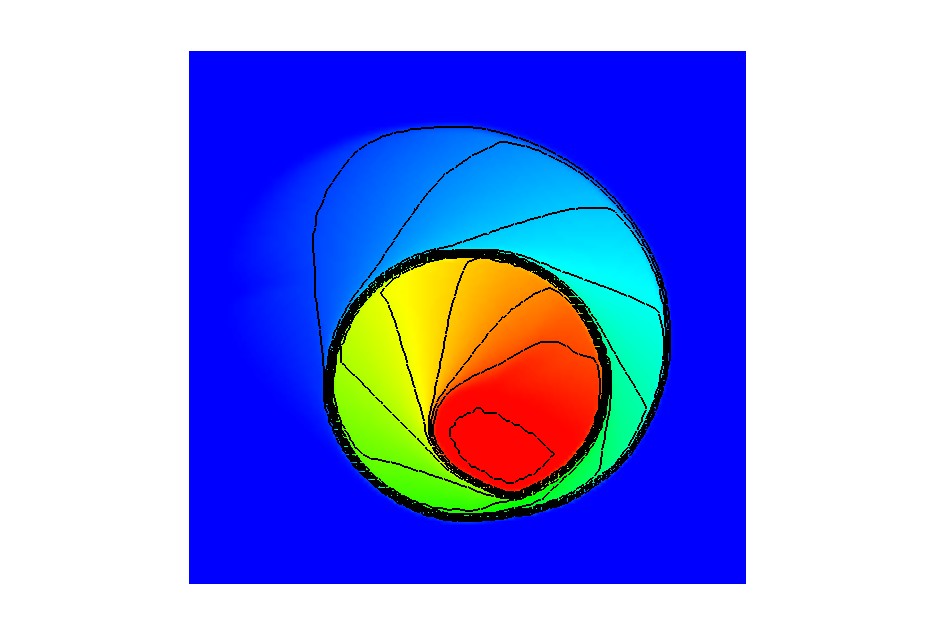}

\end{minipage}

\vskip0.35cm

\begin{minipage}[t]{0.33\textwidth}
  \centering\small ES1, $u_h\in[-0.066,12.316]$

\includegraphics[width=0.9\textwidth,trim=100 20 100 20,clip]{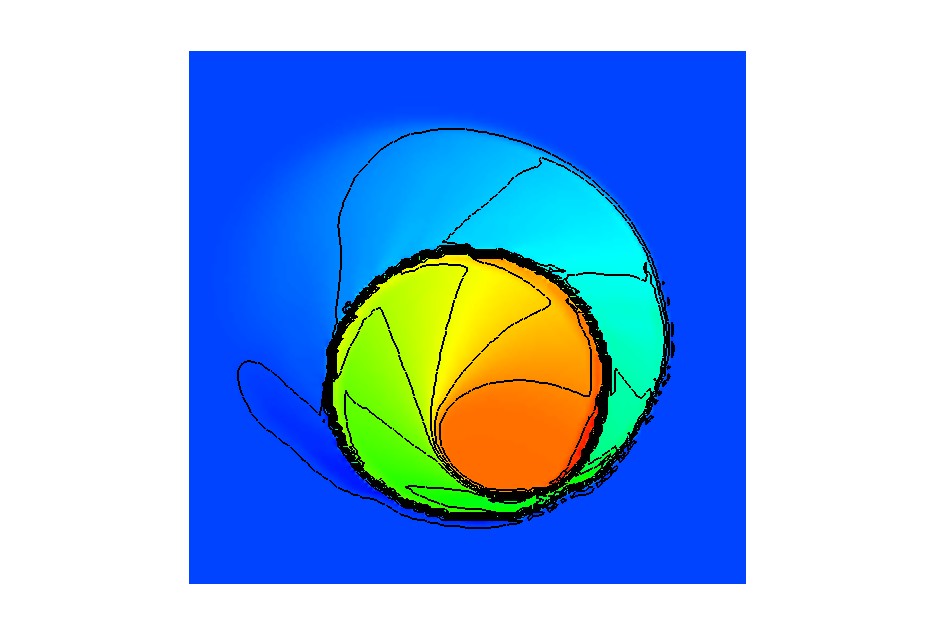}

\end{minipage}%
\begin{minipage}[t]{0.33\textwidth}

\centering\small  ES2, $u_h\in[0.037,12.241]$

\includegraphics[width=0.9\textwidth,trim=100 20 100 20,clip]{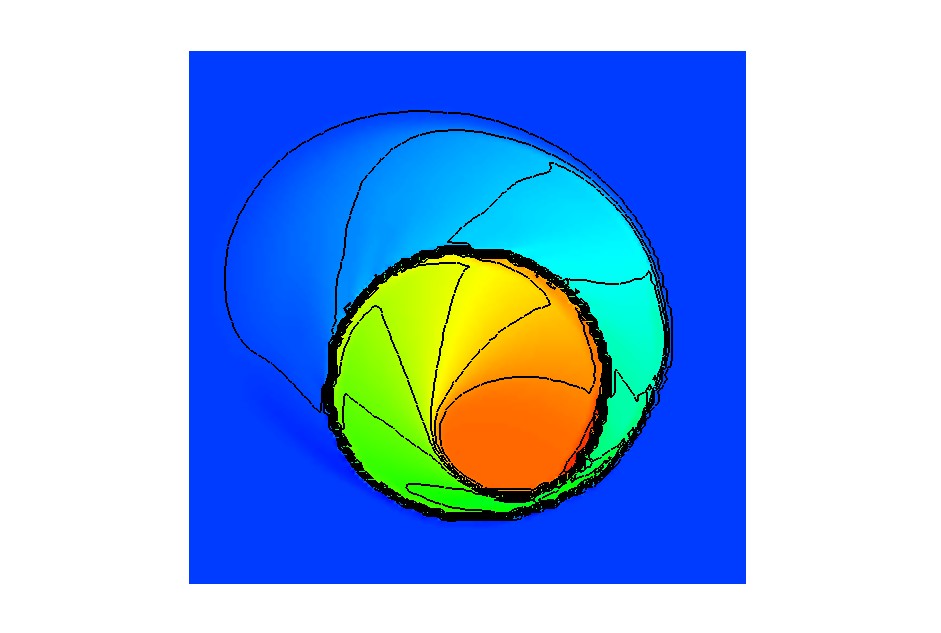}

\end{minipage}%
\begin{minipage}[t]{0.33\textwidth}

\centering\small  ES3, $u_h\in[0.213,11.504]$

\includegraphics[width=0.9\textwidth,trim=100 20 100 20,clip]{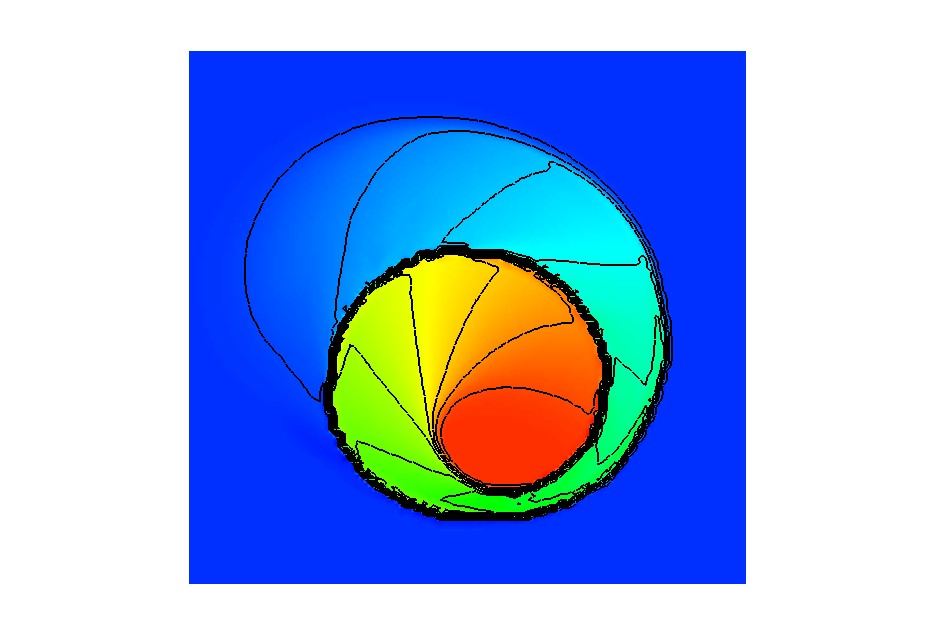}

\end{minipage}

\vskip0.35cm

\begin{minipage}[t]{0.33\textwidth}
  \centering\small  ES1S,  $u_h\in[0.785,10.996]$

\includegraphics[width=0.9\textwidth,trim=100 20 100 20,clip]{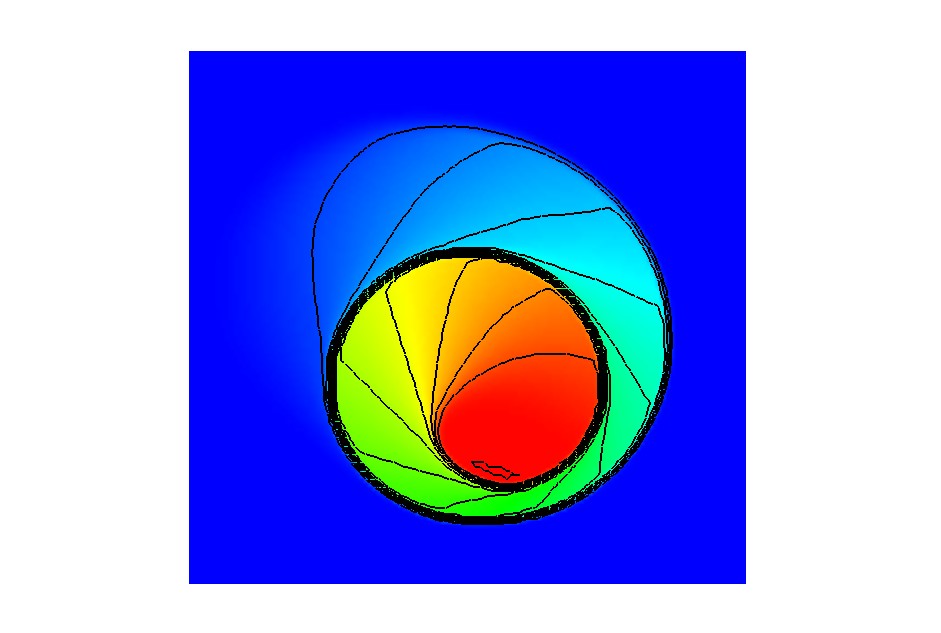}

\end{minipage}%
\begin{minipage}[t]{0.33\textwidth}

\centering\small  ES2S, $u_h\in[0.785,10.996]$

\includegraphics[width=0.9\textwidth,trim=100 20 100 20,clip]{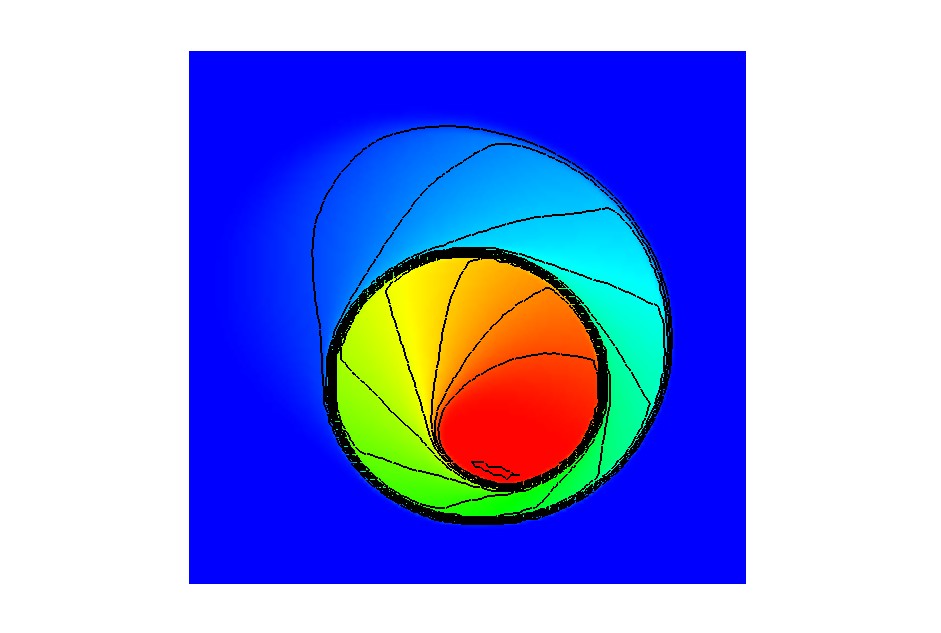}

\end{minipage}%
\begin{minipage}[t]{0.33\textwidth}

  \centering\small ES3S, $u_h\in[0.785,10.996]$

\includegraphics[width=0.9\textwidth,trim=100 20 100 20,clip]{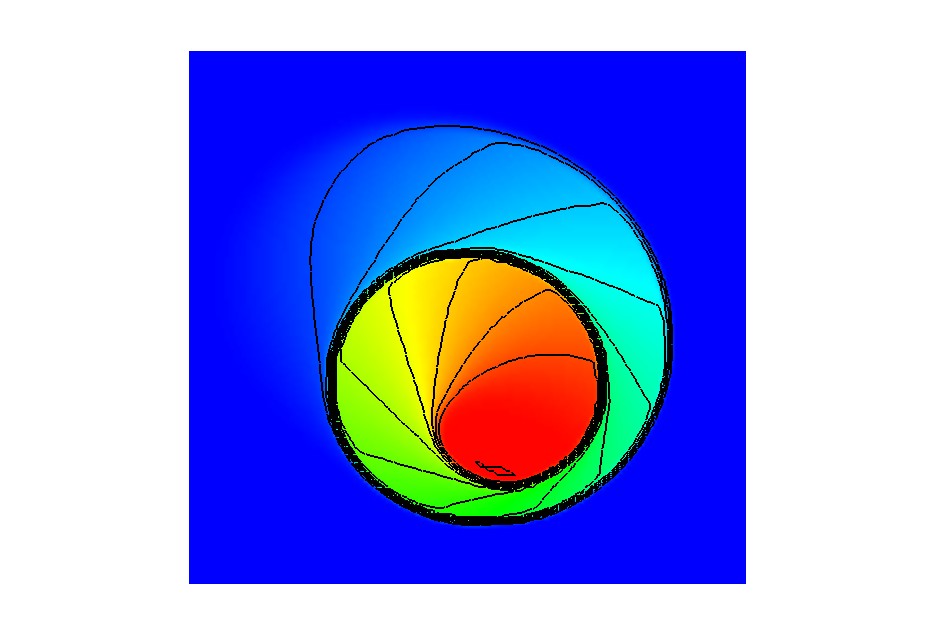}

\end{minipage}

\vskip0.35cm

\begin{minipage}[t]{0.33\textwidth}
  \centering\small  ES1FS,  $u_h\in[0.785,10.996]$

\includegraphics[width=0.9\textwidth,trim=100 20 100 20,clip]{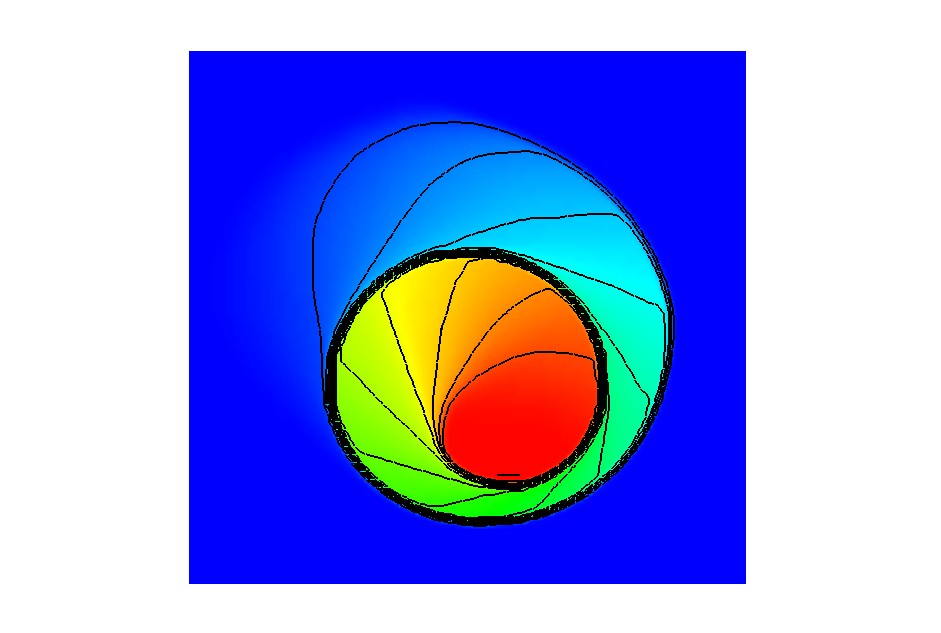}

\end{minipage}%
\begin{minipage}[t]{0.33\textwidth}

\centering\small ES2FS, $u_h\in[0.785,10.996]$

\includegraphics[width=0.9\textwidth,trim=100 20 100 20,clip]{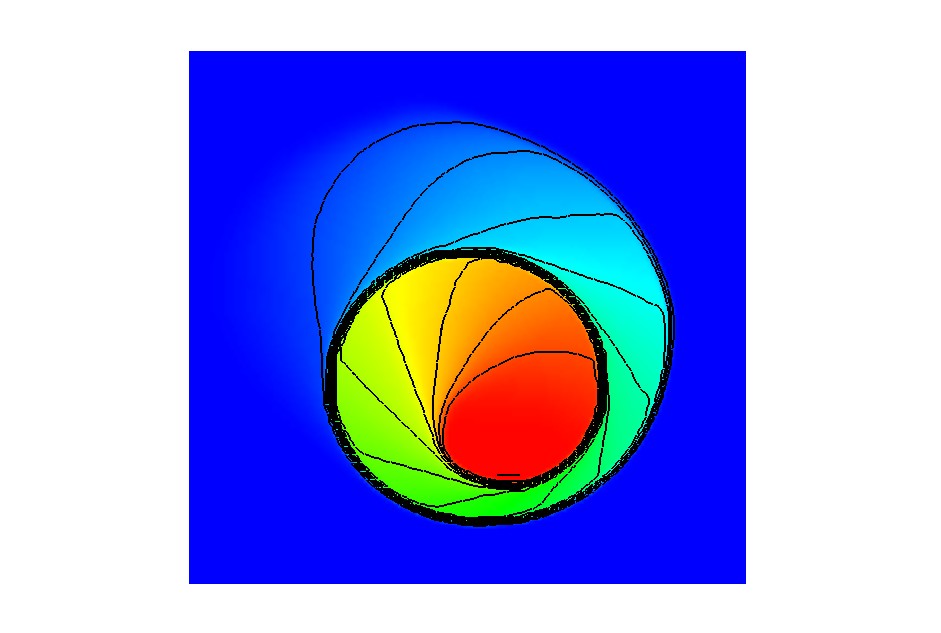}

\end{minipage}%
\begin{minipage}[t]{0.33\textwidth}

  \centering\small ES3FS, $u_h\in[0.785,10.996]$

\includegraphics[width=0.9\textwidth,trim=100 20 100 20,clip]{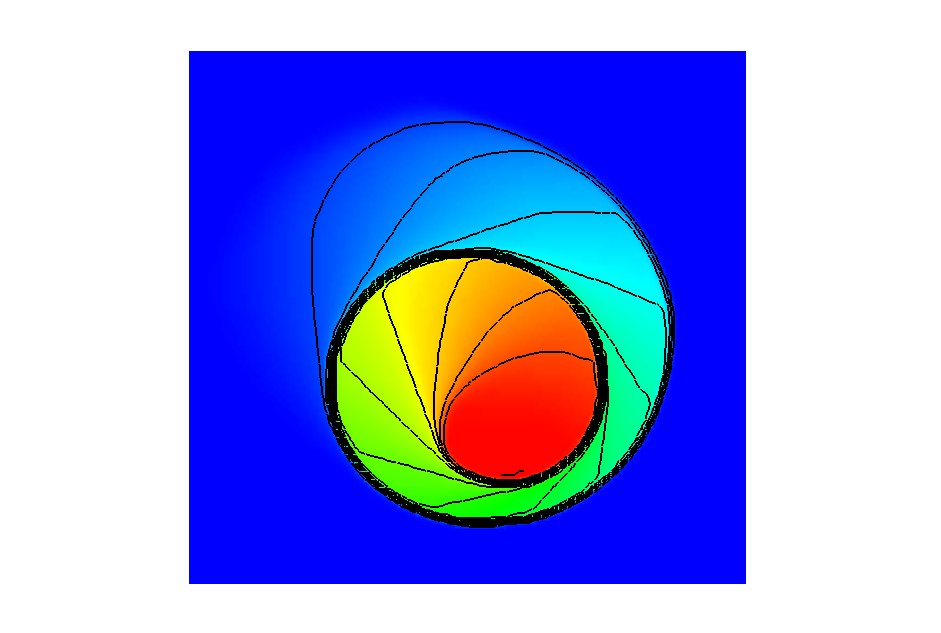}

\end{minipage}

\vskip0.25cm

\caption{KPP problem, DG solutions at $t=1.0$ calculated
  using $h=\frac{1}{128}$, $\Delta t=10^{-3}$.}

\label{KPP1}

\end{figure}

\begin{figure}[h!]
  \small
  
\begin{minipage}[t]{0.33\textwidth}
\centering DG0\vskip0.2cm

\includegraphics[width=\textwidth,trim=50 20 50 20,clip]{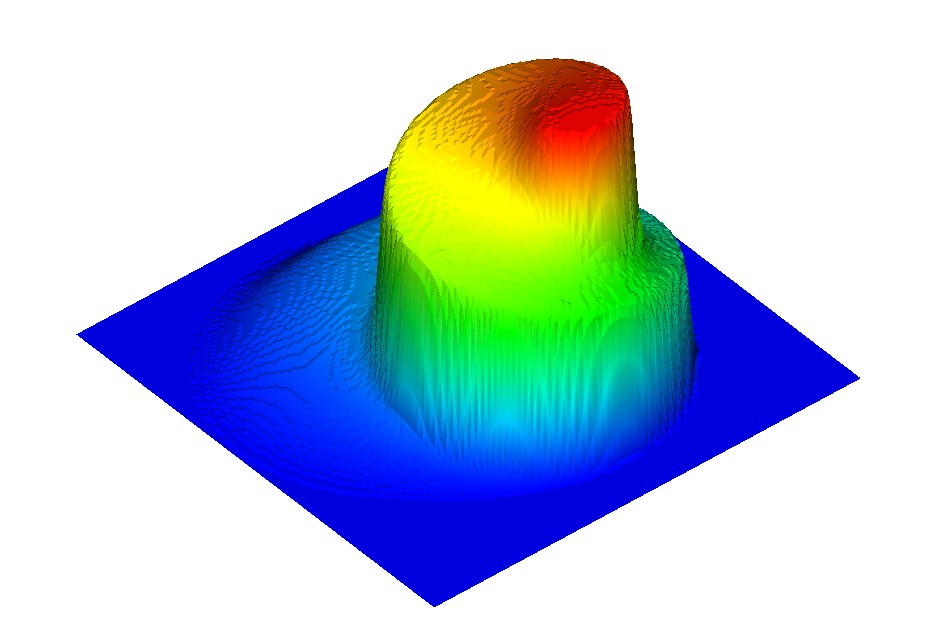}

\end{minipage}%
\begin{minipage}[t]{0.33\textwidth}

\centering DG1S\vskip0.2cm

\includegraphics[width=\textwidth,trim=50 20 50 20,clip]{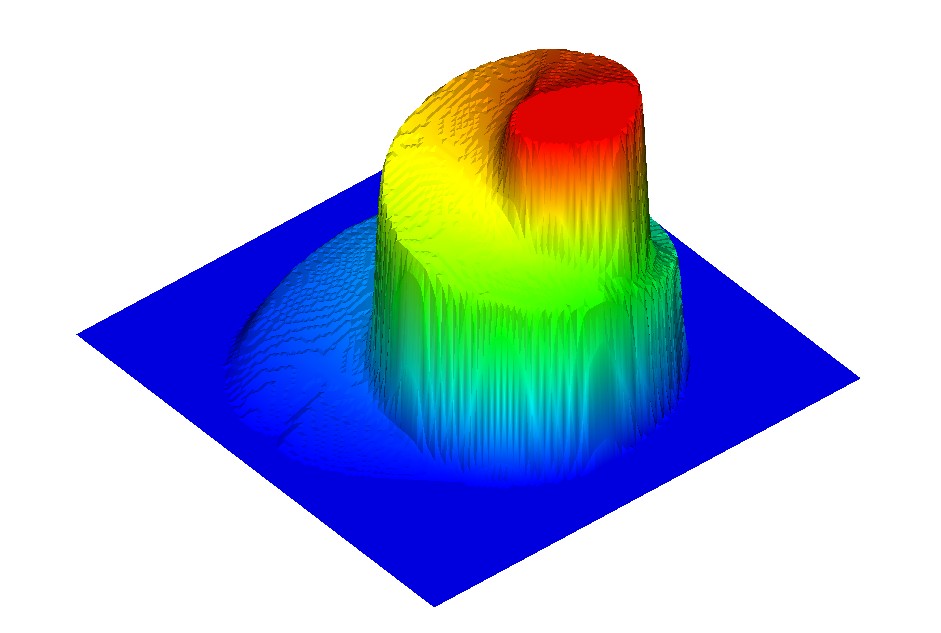}

\end{minipage}%
\begin{minipage}[t]{0.33\textwidth}

\centering ES1S\vskip0.2cm

\includegraphics[width=\textwidth,trim=50 20 50 20,clip]{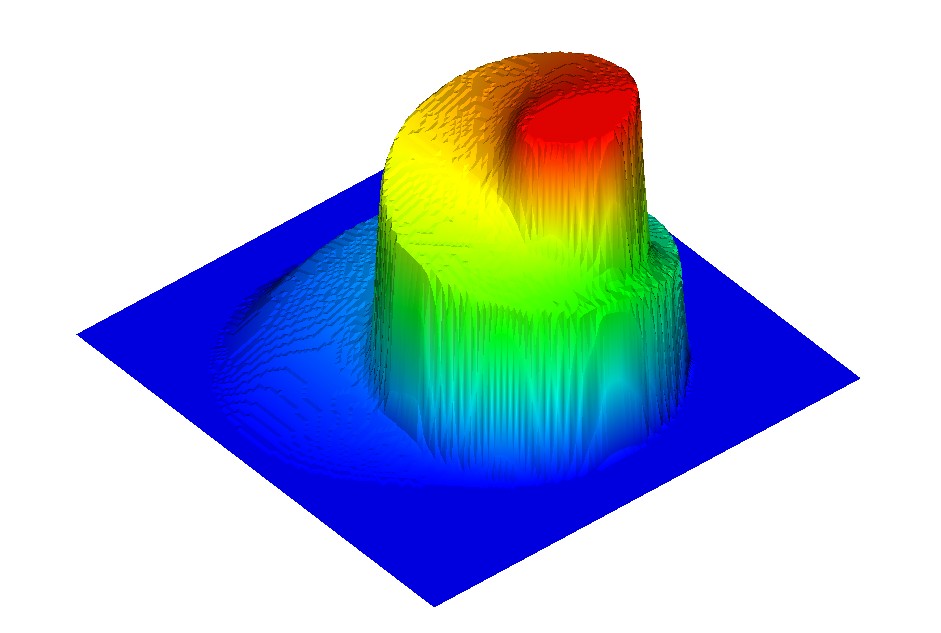}

\end{minipage}

\vskip0.25cm

\caption{KPP problem, DG solutions at $t=1.0$ calculated using
 $h=\frac{1}{128}$ and $\Delta t=10^{-3}$.}

\label{KPP2}
\vskip0.75cm
  
\begin{minipage}[t]{0.33\textwidth}
\centering DG0\vskip0.2cm

\includegraphics[width=\textwidth,trim=50 20 50 20,clip]{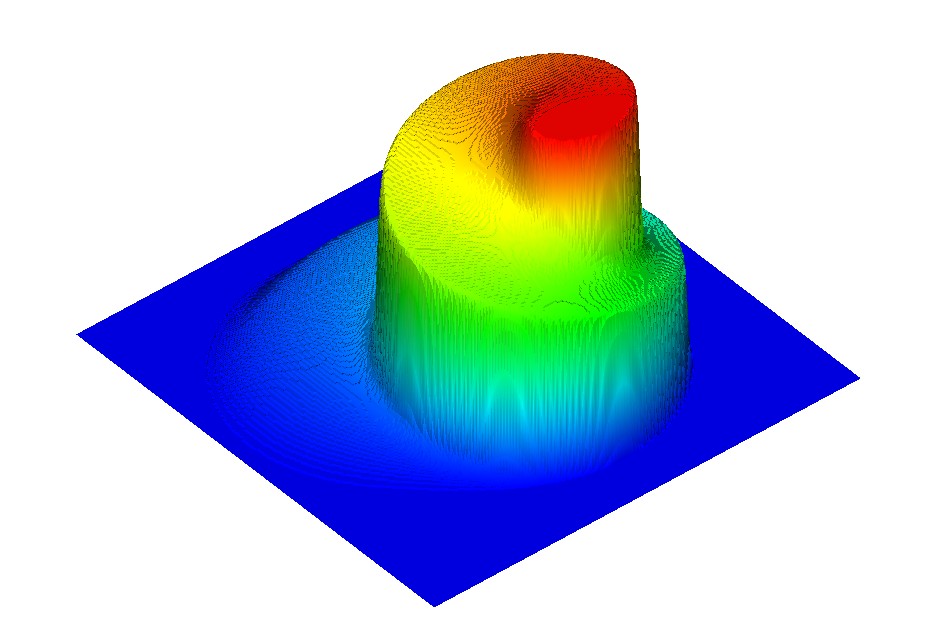}

\end{minipage}%
\begin{minipage}[t]{0.33\textwidth}

\centering DG1S\vskip0.2cm

\includegraphics[width=\textwidth,trim=50 20 50 20,clip]{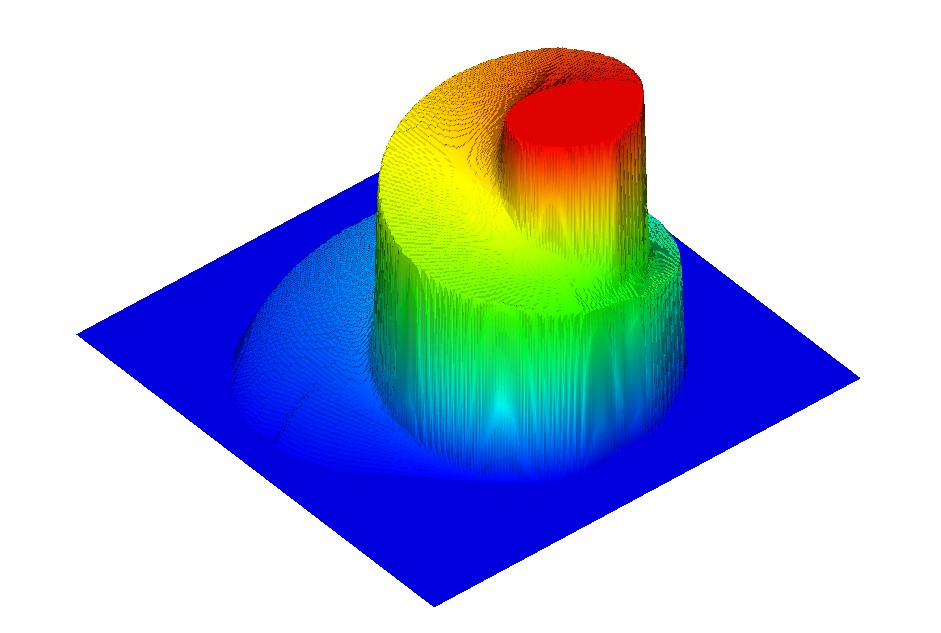}

\end{minipage}%
\begin{minipage}[t]{0.33\textwidth}

\centering ES1S\vskip0.2cm

\includegraphics[width=\textwidth,trim=50 20 50 20,clip]{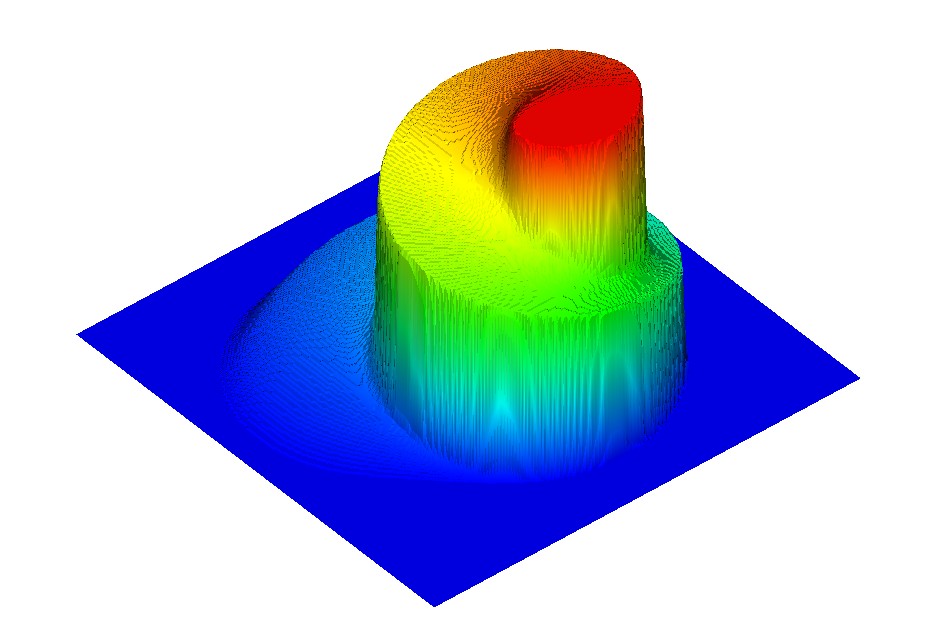}

\end{minipage}

\vskip0.25cm

\caption{KPP problem, DG solutions at $t=1.0$ calculated using 
 $h=\frac{1}{256}$ and $\Delta t=5\cdot 10^{-4}$.}

\label{KPP3}

\end{figure}

Snapshots of numerical solutions at the final time $t=1.0$ are shown
in Figs. \ref{KPP1}--\ref{KPP3}. To demonstrate the effect of entropy
stabilization and limiting, we compare the results produced by 
standard and entropy stable DG schemes without and with activation
of bound-preserving limiters (see Fig.~\ref{KPP1}). The DG0 solution is
highly dissipative but provides a correct qualitative description of
the rotating wave structure. The oscillatory DG1 solution exhibits not
only large undershoots and overshoots but also an entropy-violating
merger of two shocks. The application of the vertex-based slope
limiter in the DG1S version has the same stabilizing effect as
gradient penalization via the Rusanov dissipation term $\nu_i
\eta''(u_{i0})\bar u_{ik}$. However, the DG1 fluxes $F_{ij,k}^H$
violate condition \eqref{condES} and the semi-discrete target
scheme is not entropy stable. As a consequence of insufficient
entropy stabilization, spurious distortions are observed in the
contour lines of the slope-limited DG1S solution. Moreover, the
levels of entropy dissipation are barely enough to keep the
twisted shocks separated. This unsatisfactory state of affairs
illustrates the need for entropy stabilization and confirms 
the findings of Guermond et al. \cite{Guermond2016} who noticed
that preservation of invariant domains does not guarantee
convergence to entropy solutions.

The entropy stable schemes ES1, ES2,
ES3 preserve the thin gap between the two shocks even without
slope limiting. As the levels of entropy viscosity are increased
by decreasing the value of $Q_{ij}$ in formula \eqref{alphaES}, the
distance between the shocks increases and violations of global
bounds become less pronounced. All DG-$\mathbb{P}_1$ solutions
calculated using
bound-preserving flux and/or slope limiters look alike. They are
well-resolved and free of undershoots/overshoots. It
is not unusual that the vertex-based slope limiter produces
such solutions even if no flux limiting is performed to
constrain the local range of the cell
averages $\bar u_{j0}$ that define the bounds for 
\eqref{vblmp}. However, the
flux-limited version is generally safer because the
validity of local maximum principles is guaranteed by
Theorem~\ref{thm2}.

The snapshots shown in Fig. \ref{KPP2} visualize the corresponding
diagrams of  Fig. \ref{KPP1} as surface plots to better illustrate
the capability of DG0, DG1S, and ES1S to capture the rotating
wave structure on the mesh with spacing $h=\frac{1}{128}$.
The results obtained with the
three methods on a finer mesh ($h=\frac{1}{256}$)
are displayed in Fig. \ref{KPP3}. 
The fine-mesh DG0 and ES1S solutions illustrate the correct
shock behavior. The DG1S solutions indicate that the use of
slope limiting has a strong stabilizing effect but may fail
to prevent entropy-violating behavior that does not cause
violations of local bounds in \eqref{vblmp}.

\subsection{Buckley-Leverett equation}
In the second numerical experiment, we consider the two-dimensional
Buckley-Leverett equation \cite{christov2008new,entropyCG}. The nonconvex
flux function of the nonlinear conservation law to be solved
is  
\beq
\mathbf{f}(u)=
\frac{u^2}{u^2+(1-u)^2}
(1,1-5(1-u)^2).
\eeq
The computational domain is $\Omega_h=(-1.5,1.5)^2$. The
piecewise-constant initial condition is given by
\beq
u_0(x,y)=\begin{cases}
1 & \mbox{ if } x^2+y^2<0.5, \\
0 & \mbox{ otherwise}.
\end{cases}
\eeq
Similarly to the KPP problem, the solution exhibits a rotating wave structure. For entropy stabilization purposes, we use $\eta(u)=\frac{u^2}2$ and the corresponding entropy flux
$\mathbf{q}(u)=(q_x(u),q_y(u))$, where
\begin{align}
q_x&=\frac14\left[
\frac{2(u-1)}{2u^2-2u+1}-\log(2u^2-2u+1)\right],\\
q_y&=\frac{1}{12}\left[
-20u^3+15u^2-\frac{9u+6}{2u^2-2u+1}-3\log(2u^2-2u+1)
-15\tan^{-1}(1-2u)\right].
\end{align}
An upper bound for the fastest wave speed can be found in \cite{christov2008new}. We overestimate it by using $\lambda_{ij}=3.4$.

\begin{figure}[h!]
  
\begin{minipage}[t]{0.33\textwidth}
  \centering\small DG0  $u_h\in[0.0,0.959]$

\includegraphics[width=0.9\textwidth,trim=100 20 100 20,clip]{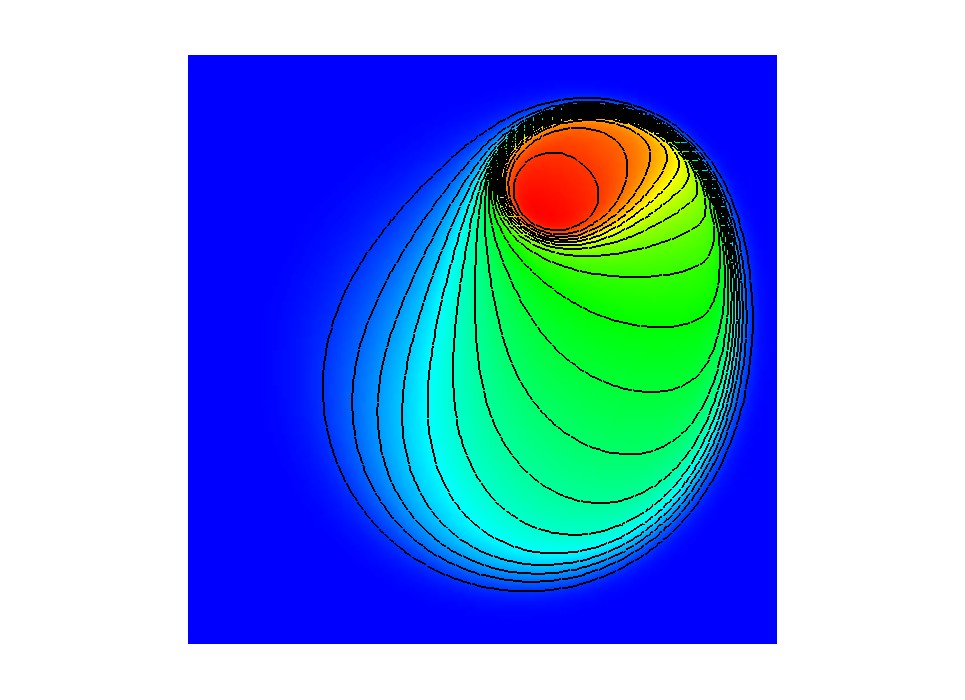}

\end{minipage}%
\begin{minipage}[t]{0.33\textwidth}

\centering\small  DG1, $u_h\in[-0.268,1.268]$

\includegraphics[width=0.9\textwidth,trim=100 20 100 20,clip]{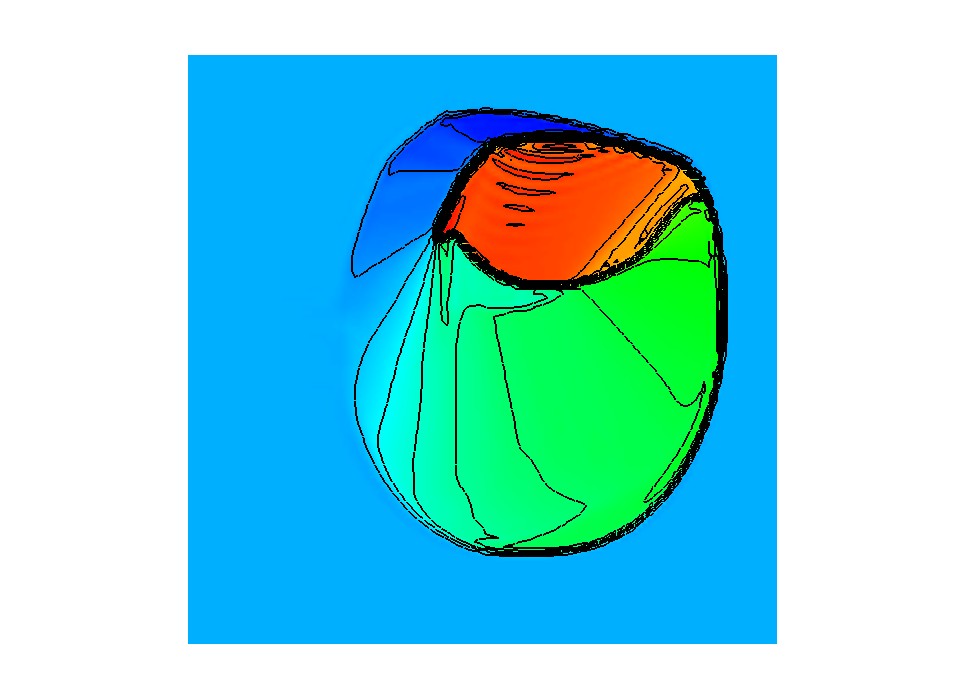}

\end{minipage}%
\begin{minipage}[t]{0.33\textwidth}

  \centering\small  DG1S, $u_h\in[-0.354,1.451]$

\includegraphics[width=0.9\textwidth,trim=100 20 100 20,clip]{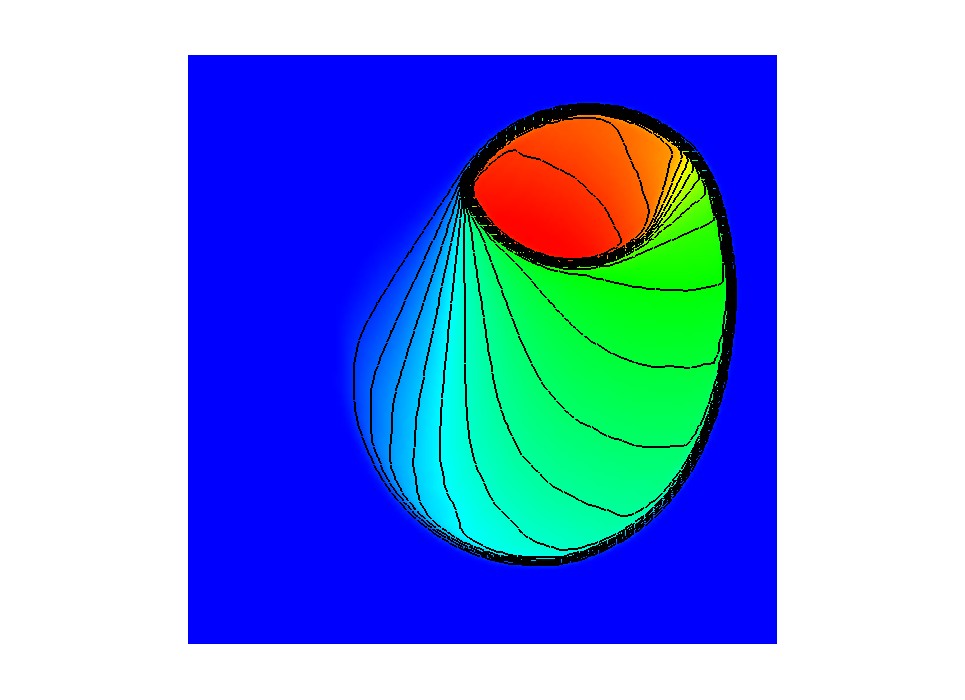}

\end{minipage}

\vskip0.35cm

\begin{minipage}[t]{0.33\textwidth}
  \centering\small ES1, $u_h\in[-0.046,1.090]$

\includegraphics[width=0.9\textwidth,trim=100 20 100 20,clip]{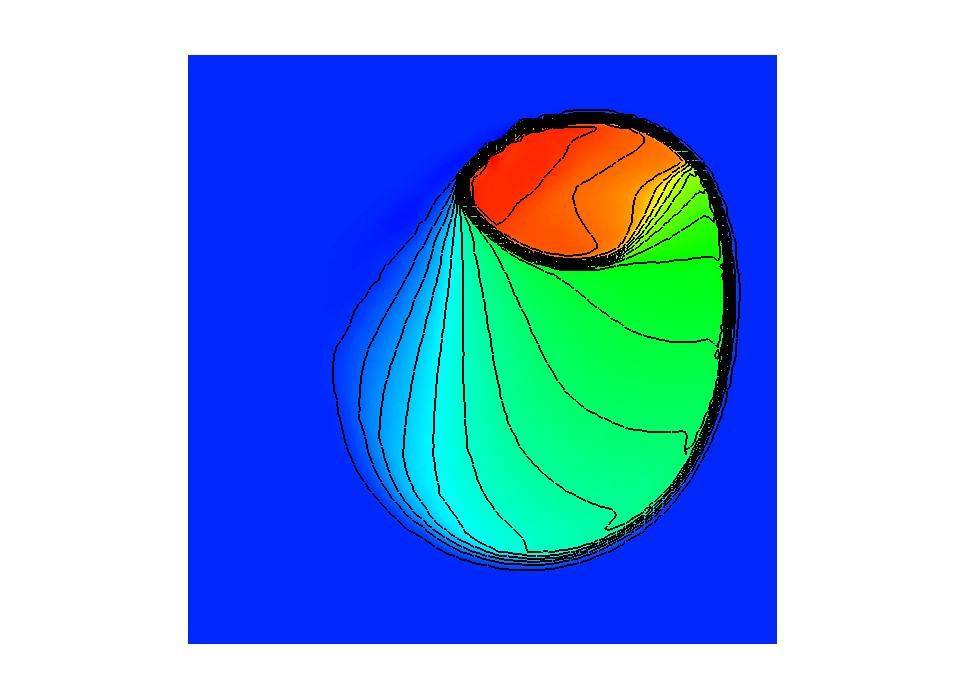}

\end{minipage}%
\begin{minipage}[t]{0.33\textwidth}

\centering\small  ES2, $u_h\in[-0.083,1.149]$

\includegraphics[width=0.9\textwidth,trim=100 20 100 20,clip]{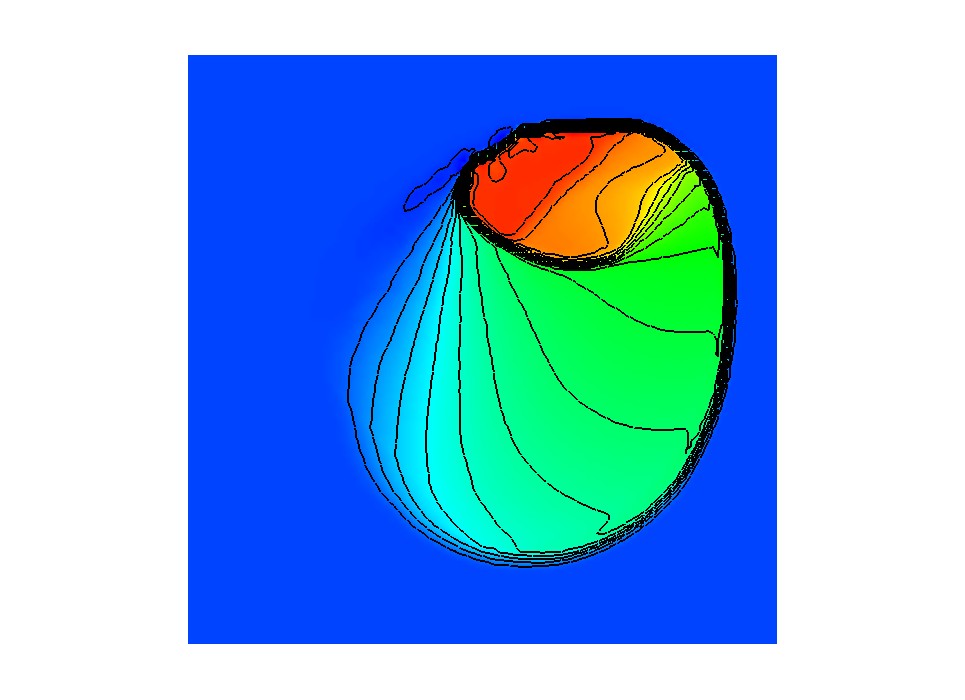}

\end{minipage}%
\begin{minipage}[t]{0.33\textwidth}

\centering\small  ES3, $u_h\in[-0.015,1.140]$

\includegraphics[width=0.9\textwidth,trim=100 20 100 20,clip]{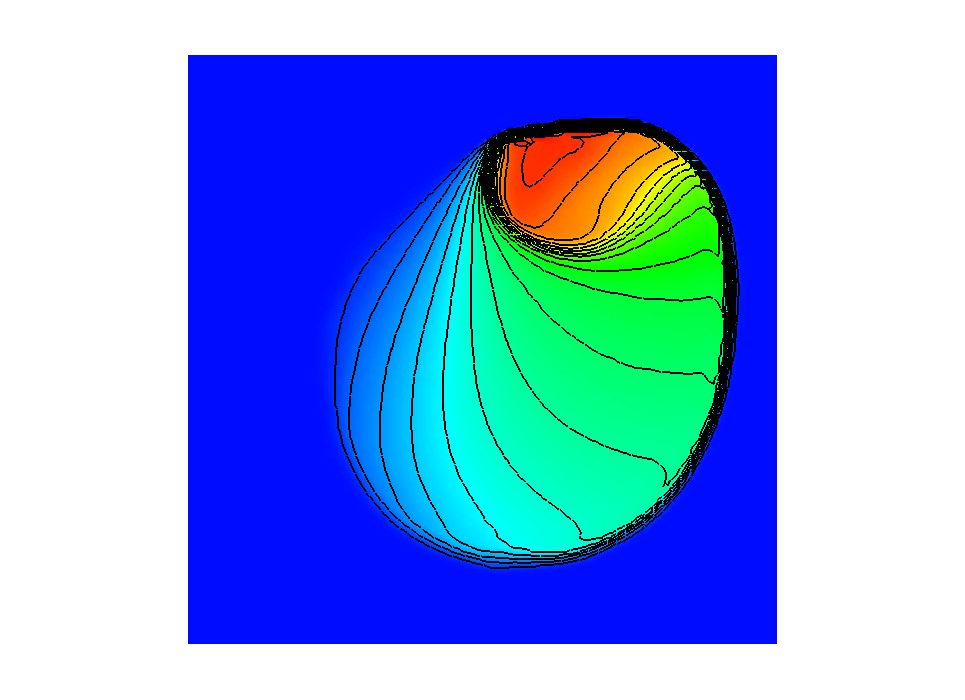}

\end{minipage}

\vskip0.35cm

\begin{minipage}[t]{0.33\textwidth}
  \centering\small  ES1S,  $u_h\in[0.0,0.997]$

\includegraphics[width=0.9\textwidth,trim=100 20 100 20,clip]{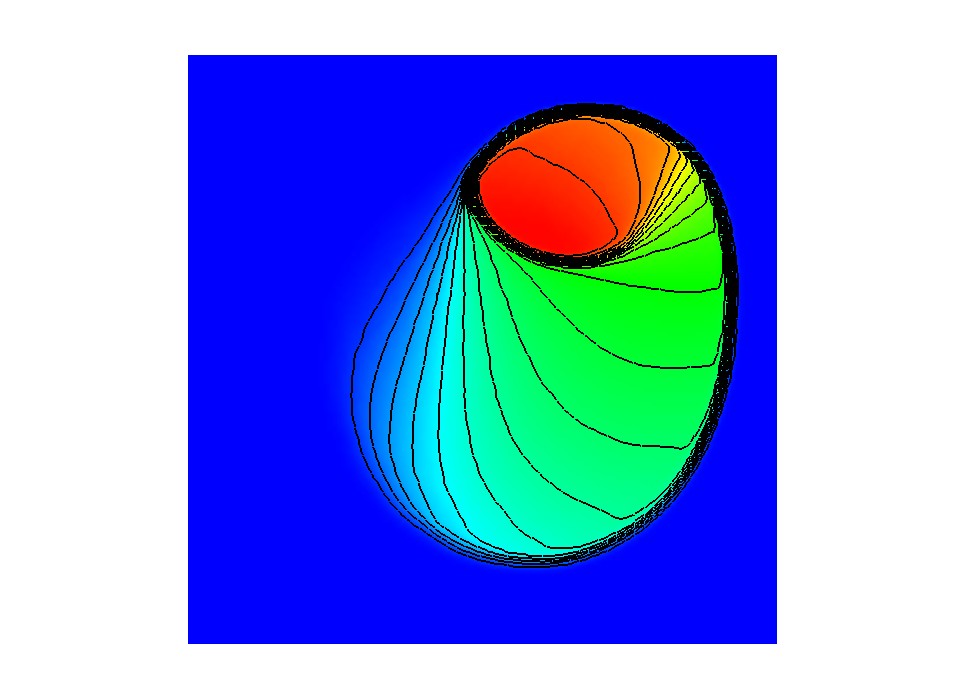}

\end{minipage}%
\begin{minipage}[t]{0.33\textwidth}

\centering\small  ES2S, $u_h\in[0.0,0.997]$

\includegraphics[width=0.9\textwidth,trim=100 20 100 20,clip]{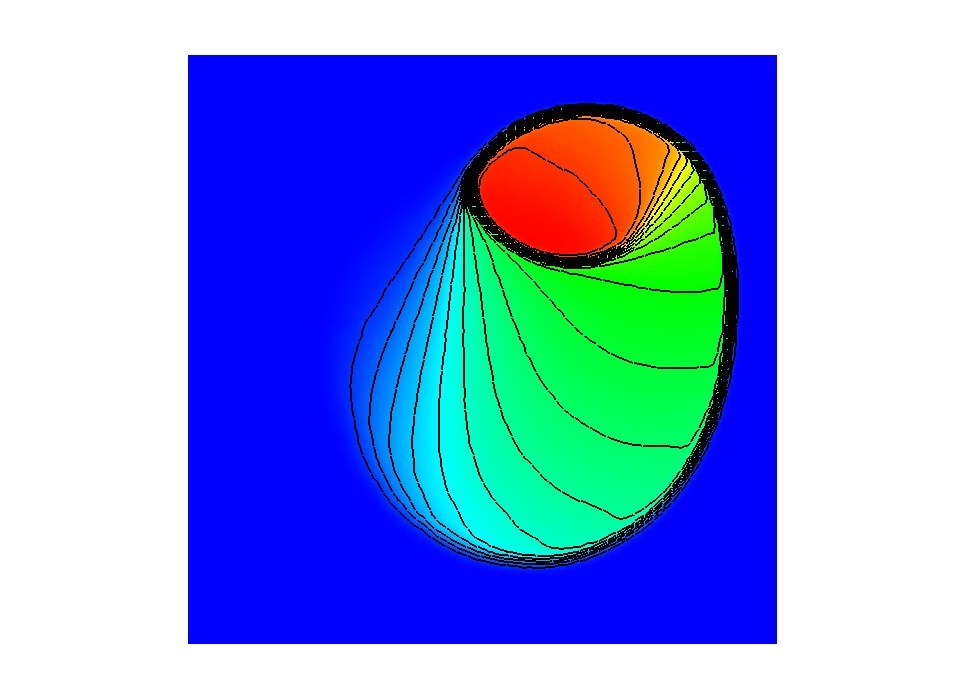}

\end{minipage}%
\begin{minipage}[t]{0.33\textwidth}

  \centering\small ES3S, $u_h\in[0.0,0.990]$

\includegraphics[width=0.9\textwidth,trim=100 20 100 20,clip]{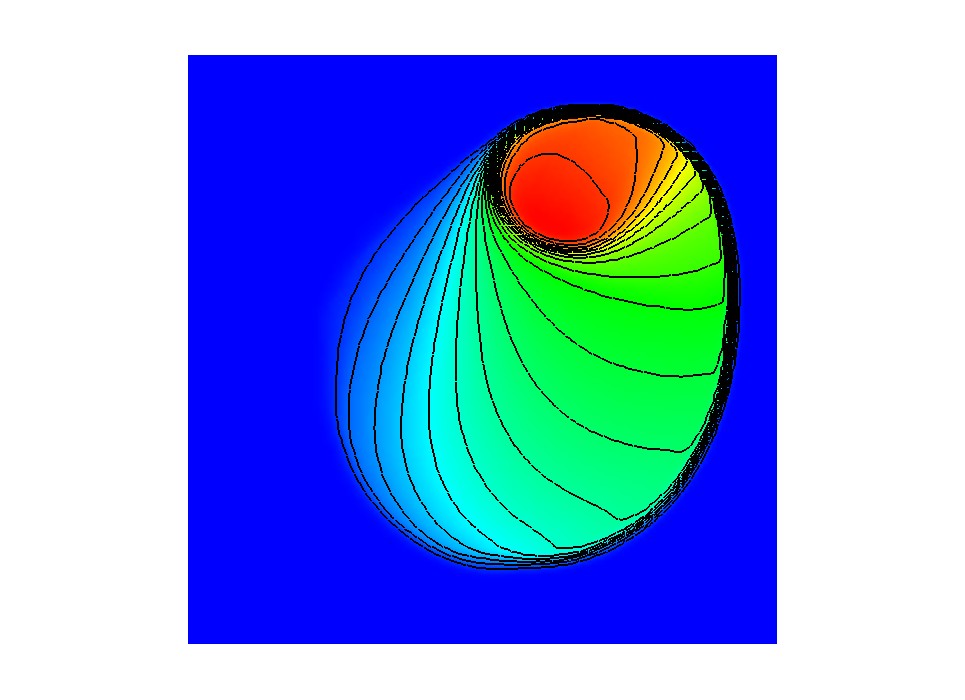}

\end{minipage}

\vskip0.25cm

\caption{Buckley-Leverett problem, DG solutions at $t=0.5$ calculated
  using $h=\frac{1}{128}$, $\Delta t=10^{-3}$.}

\label{BL}

\end{figure}

In Figure \ref{BL}, we show the numerical results at
the final time $t=0.5$  obtained
using a uniform mesh with spacing
$h=\frac{1}{128}$. The qualitative behavior of the
DG solutions is similar to that
for the more challenging KPP problem. The 
DG0 approximation is bound-preserving w.r.t.  
$\mathcal G=[0,1]$ but the levels of numerical diffusion
are too high. The
DG1 scheme produces an oscillatory solution and gives rise to
strong violations of the global bounds.
The entropy stable approximations ES1, ES2, ES3
exhibit smaller undershoots and overshoots. The application of
the vertex-based slope limiter eliminates them completely. The
results obtained with the bound-preserving flux limiter look
similar and are not presented here.

\section{Conclusions}
\label{sec:end}

The main result of the presented work is the development of
a methodology that pieces together individual components of
property-preserving DG schemes for nonlinear hyperbolic
problems. We have shown that
a carefully tuned combination of flux limiting, entropy
stabilization, and slope limiting makes it possible to
satisfy all relevant inequality constraints by blending
a low-order LLF approximation and a high-order target.
Although the proposed framework was introduced in the
context of scalar nonlinear conservation laws, its
extension to systems appears to be straightforward.
In the case of $m>1$ conserved quantities, a generalized
version of the entropy condition \eqref{condES} and
preservation of invariant domains can be readily enforced using
\eqref{fluxlim} with a scalar correction factor $\alpha_{ij}$.
Additionally, the components of the high-order target flux
$F_{ij,0}^H$ can be pre-constrained to satisfy local maximum
principles for certain quantities of interest (see
\cite{convex} for details). Last but not least, the time
integration procedure can be redesigned to be not only SSP
but also entropy stability preserving \cite{relax}. It
is hoped that the findings of this paper will provide useful
insights and tools for such research endeavors.

\medskip
\paragraph{\bf Acknowledgments}

This research was supported by the German Research Association (DFG) under grant KU 1530/23-1.  The author would like to thank Manuel Quezada de Luna (KAUST) and Hennes Hajduk (TU Dortmund University) for inspiring discussions of the presented methodology.


\end{document}